\documentclass[preprint]{elsarticle}

\usepackage{csl}
\usepackage{epsfig, graphicx}
\usepackage{makeidx}
\usepackage{multicol,multirow}
\usepackage{amsmath, amsfonts, amssymb, bbm, bm, mathtools, dsfont}
\usepackage{ulem}
\usepackage{crop}
\usepackage{xcolor,color}
\usepackage{todonotes}
\usepackage{hyperref,url}
\usepackage{cleveref,scalerel}
\usepackage{xspace,ifthen,mleftright,scalerel}
\usepackage{algorithmic}

\usepackage{tikz}
\usepackage{tikz-qtree}
\usepackage{subfiles}
\usepackage{stmaryrd}
\usepackage{pdflscape}
\usepackage{breqn}
\usepackage{rotating}

\usepackage{array}
\newcolumntype{C}[1]{>{\centering\arraybackslash}m{#1}}  
\newcolumntype{L}[1]{>{\raggedleft\arraybackslash}m{#1}} 
\newcolumntype{R}[1]{>{\raggedleft\arraybackslash}m{#1}} 

\usepackage{pgfplots}
\pgfplotsset{compat=1.14}
\pgfplotsset{every axis/.append style={
    legend style={at={(0.5,1.05)},anchor=south},
    legend style={/tikz/every even column/.append style={column sep=1cm}},
    legend style={font=\scriptsize},
    legend columns=3
}}

\usepackage{cite}

\usepackage{xparse}

\newcommand{\slow}{{\scriptstyle \textsc{s}}}
\newcommand{\fast}{{\scriptstyle \textsc{f}}}

\newcommand{\crl}[1]{\raisebox{0.5pt}{\scaleto{\{}{5pt}}\hspace*{-0.25pt} #1\hspace*{-0.25pt}\raisebox{0.5pt}{\scaleto{\}}{5pt}}}

\newcommand{\F}[1][]{%
   \ifthenelse{ \equal{#1}{} }
      {^{\crl{\fast}}}
      {^{\crl{\fast,\scaleto{#1}{4pt}}}}
}
\newcommand{\FF}[1][]{%
   \ifthenelse{ \equal{#1}{} }
      {^{\crl{\fast,\fast}}}
      {^{\crl{\fast,\fast,\scaleto{#1}{4pt}}}}
}
\newcommand{\FS}[1][]{%
   \ifthenelse{ \equal{#1}{} }
      {^{\crl{\fast,\slow}}}
      {^{\crl{\fast,\slow,\scaleto{#1}{4pt}}}}
}
\newcommand{\SF}[1][]{%
   \ifthenelse{ \equal{#1}{} }
      {^{\crl{\slow,\fast}}}
      {^{\crl{\slow,\fast,\scaleto{#1}{4pt}}}}
}

\newcommand{\nvar}{d}

\newcommand{\fun}{\mathbf{f}}

\newcommand{\y}{\mathbf{y}}

\newcommand{\yf}[1][]{%
   \ifthenelse{ \equal{#1}{} }
      {{\mathbf{y}_\fast}}
      {{\mathbf{y}_{\fast,#1}}}
}

\newcommand{\ys}[1][]{%
   \ifthenelse{ \equal{#1}{} }
      {{\mathbf{y}_\slow}}
      {{\mathbf{y}_{\slow,#1}}}
}

\newcommand{\Id}{\mathbf{I}}
\newcommand{\Idvar}{{\Id_{\nvar}}} 
\newcommand{\Idstage}{{\Id_{s}}} 
\newcommand{\one}{1\hspace{-0,9ex}1} 

\def\Re{\mathds{R}}

\newcommand{\diag}{\operatorname{diag}}

\newcommand{\tree}{\mathfrak{t}}
\newcommand{\utree}{\mathfrak{u}}

\newcommand{\msquare}[1]{\scalebox{0.5}{\framebox{$#1$}}}
\newcommand{\mround}[1]{\scalebox{0.7}{\textcircled{\scalebox{0.7}{$#1$}}}}

\newtheorem{remark}{Remark}[section]

\newtheorem{definition}{Definition}[section]
\newtheorem{lemma}{Lemma}[section]
\newtheorem{theorem}{Theorem}[section]
\newtheorem{proof}{Proof}[section]

\newcommand{\A}{{\mathbf{A}}}
\renewcommand{\b}{{\mathbf{b}}}

\renewcommand{\c}{{\mathbf{c}}}

\newcommand{\nparts}{{\rm N}}
\newcommand{\NT}{\mathds{T}_{\nparts}}
\newcommand{\NTW}{\mathds{T\!W}_{\nparts}}
\newcommand{\DAT}{\mathds{D\!A\!T}}
\newcommand{\Times}{%
  \resizebox{!}{1.5\fontcharht\font`0}{$\mkern-2mu\times\mkern-2mu$}%
}

\newenvironment{butchertableau}[2][1.3]
	{\def\arraystretch{#1}\array{#2}}
	{\endarray}

\newenvironment{butchermatrix}[1][1.3]
	{\def\arraystretch{#1}\begin{bmatrix}}
	{\end{bmatrix}}

\newcommand{\G}{\mathbf{G}}
\renewcommand{\b}{\mathbf{b}}
\renewcommand{\c}{\mathbf{c}}

\newcommand{\Lb}{\mathbf{L}}
\newcommand{\Jac}{\mathbf{J}}

\newcommand{\impl}{\textsc{i}}
\newcommand{\expl}{\textsc{e}}
\newcommand{\I}{^{\{\impl\}}}
\newcommand{\E}{^{\{\expl\}}}
\newcommand{\IE}{^{\{\impl,\expl\}}}
\newcommand{\II}{^{\{\impl,\impl\}}}
\newcommand{\EE}{^{\{\expl,\expl\}}}
\newcommand{\EI}{^{\{\expl,\impl\}}}

\newcommand{\alg}{{\rm z}}
\newcommand{\diff}{{\rm x}}
\newcommand{\Alg}{^{\{\alg\}}}
\newcommand{\Diff}{^{\{\diff\}}}
\newcommand{\AlgDiff}{^{\{\alg,\diff\}}}
\newcommand{\AlgAlg}{^{\{\alg,\alg\}}}
\newcommand{\DiffDiff}{^{\{\diff,\diff\}}}
\newcommand{\DiffAlg}{^{\{\diff,\alg\}}}

\newcommand{\sfrac}[2]{\mbox{\footnotesize$\displaystyle\frac{#1}{#2}$}}

\newcommand{\gun}{\mathbf{g}}

\newcommand{\xx}{\mathbf{x}}
 \newcommand{\yy}{\mathbf{y}}
 \newcommand{\zz}{\mathbf{z}}

\makeatletter
\newcommand{\kron}[1]{{\,\mathrlap{\otimes}{\mathrlap{\hspace{0.275em}\tikz{\path[draw=white,fill=white] (0.025em,0.025em) rectangle (0.22em,0.38em);}}}\scalebox{0.4}{\raisebox{3.5pt}{\hspace{0.7em}#1}} \hspace{0.5em} }}
\makeatother

\newif\ifreport
\reporttrue

\begin{document}
	
	\ifreport
	\csltitle{Linearly implicit GARK schemes}
	\cslauthor{Adrian Sandu, Michael G\"unther, and Steven Roberts}
	\cslemail{sandu@cs.vt.edu, guenther@uni-wuppertal.de, steven94@vt.edu}
	\cslreportnumber{9}
	\cslyear{20}
	\csltitlepage
	\fi
	
\begin{frontmatter}

\title{Linearly implicit GARK schemes\tnoteref{t1}}
\tnotetext[t1]{The work of Sandu and Roberts was supported by awards NSF CCF--1613905. NSF ACI--1709727, NSF CDS\&E-MSS--1953113, and by the Computational Science Laboratory at Virginia Tech. The work of G\"unther was supported  by the European Union's Horizon 2020 research and innovation programme under the Marie Sklodowska-Curie Grant Agreement No.~765374, ROMSOC. 
}

\author[1]{Adrian Sandu}
\ead{sandu@cs.vt.edu}
\author[2]{Michael G\"unther\corref{cor1}}
\ead{guenther@uni-wuppertal.de}
\author[3]{Steven Roberts}
\ead{steven94@vt.edu}
\cortext[cor1]{Corresponding author}
\address[1]{Computational Science Laboratory, Department of  Computer Science, Virginia Tech, Blacksburg, VA 24061}
\address[2]{IMACM, School of Mathematics and Natural Sciences, Department of Mathematics and Informatics, Bergische Universit\"at Wuppertal, 42097 Wuppertal, Germany}
\address[3]{Computational Science Laboratory, Department of  Computer Science, Virginia Tech, Blacksburg, VA 24061}

\begin{abstract}
Systems driven by multiple physical processes are central to many areas of science and engineering.  Time discretization of multiphysics systems is challenging, since different processes have different levels of stiffness and characteristic time scales.  The multimethod approach discretizes each physical process with an appropriate numerical method; the methods are coupled appropriately such that the overall solution has the desired accuracy and stability properties. The authors developed the general-structure additive Runge--Kutta (GARK) framework, which constructs multimethods based on Runge--Kutta schemes. 

This paper constructs the new GARK-ROS/GARK-ROW families of multimethods based on linearly implicit Rosenbrock/Rosenbrock-W schemes. For ordinary differential equation models, we develop a general order condition theory for linearly implicit methods with any number of partitions, using exact or approximate Jacobians.  We generalize the order condition theory to two-way partitioned index-1 differential-algebraic equations. Applications of the framework include decoupled linearly implicit, linearly implicit/explicit, and linearly implicit/implicit methods. Practical GARK-ROS and GARK-ROW schemes of order up to four are constructed.
\end{abstract}

\begin{keyword}
Multiphysics systems \sep GARK methods \sep linear implicitness
\MSC 65L05, 65L06, 65L07, 65L20.
\end{keyword}

\end{frontmatter}

\section{Introduction}

We are concerned with the numerical solution of differential equations arising in the simulation of multiphysics systems. Such equations are of great practical importance as they model diverse phenomena that appear in mechanical and chemical engineering, aeronautics, astrophysics, plasma physics, meteorology and oceanography, finance, environmental sciences, and urban modeling.
A general representation of multiphysics dynamical systems has the form:
\begin{equation}
\label{eqn:additively-partitioned-ode}
\frac{d\yy}{dt} = \fun(\yy)  = \sum_{m=1}^\nparts \fun^{\{m\}}(\yy),  ~~ t_0 \leq t \leq t_F, ~~ \yy(t_0) = \yy_0 \in \Re^\nvar,
\end{equation}
where \cref{eqn:additively-partitioned-ode} is driven by multiple physical processes $\fun^{\{m\}} : \Re^\nvar \rightarrow \Re^\nvar$ with different dynamical characteristics, and acting simultaneously.

Time discretization of complex systems \cref{eqn:additively-partitioned-ode} is challenging, since different processes have different levels of stiffness and characteristic time scales. Explicit schemes \citep{Hairer_book_I} advance the solution using only information from previous steps at a low computational cost per-timestep; however,  in addition to step size limitations   due to stability considerations, explicit timesteps can be only as large as the fastest time scale in the system. Implicit schemes that advance solutions using past and future information \citep{Hairer_book_II} remove the stability restrictions on timestep size; however their computational cost per-timestep is large, as they solve one or more systems of nonlinear equations. 
Stiffness in any individual process requires the use of an implicit solver for the entire multiphysics system \cref{eqn:additively-partitioned-ode}.

Linearly implicit methods seek to preserve the good stability properties of implicit schemes, but avoid solving large nonlinear systems of equations; instead, they only require solutions of linear systems at each step. In his seminal 1963 paper \citep{Rosenbrock_1963} Rosenbrock proposed linearly implicit Runge--Kutta type methods. 
%
%
%
An $s$-stage Rosenbrock  method solves the autonomous system \cref{eqn:additively-partitioned-ode} in its aggregated form (i.e., treating all individual components in the same way) as follows \citep[Section IV.7]{Hairer_book_II}
\begin{subequations}
\label{eqn:ROS}
\begin{align}
\label{eqn:ROS-stage}
k_i &= h \,\fun\mleft( \yy_{n}+\displaystyle\sum_{j=1}^{i-1}\alpha_{i,j} \,k_j \mright)  
 + h \,\Jac_n\, \displaystyle\sum_{j=1}^{i} \gamma_{i,j}\, k_j  	
		\,, ~~  i=1,\dots,s, \\
\label{eqn:ROS-solution}
    		\yy_{n+1} &= \yy_{n} + \displaystyle\sum_{i=1}^s b_i\, k_i,
\end{align}
\end{subequations}
where the matrix $\Jac_{n} \coloneqq \fun_\yy(\yy_{n}) \in \Re^{\nvar \times \nvar}$ is the Jacobian of the aggregated right hand side function \cref{eqn:additively-partitioned-ode}. Each stage vector $k_i$ is the solution of a linear system with matrix $\Idvar - h\,\gamma_{i,i}\,\Jac_n$, and if $\gamma_{i,i} = \gamma$ for all $i$ then the same LU factorization can be reused for all stages.
We consider the following matrices of method coefficients:
\begin{equation}
\mathbf{b} =  [ b_{i} ]_{1 \le i \le s}, \quad
\boldsymbol{\alpha} = [ \alpha_{i,j} ]_{1 \le i,j \le s}, \quad
\boldsymbol{\gamma} = [ \gamma_{i,j} ]_{1 \le i,j \le s}, \quad
\boldsymbol{\beta} = \boldsymbol{\alpha} + \boldsymbol{\gamma},
\end{equation}
where in \cref{eqn:ROS-stage} $\boldsymbol{\alpha}$ is strictly lower triangular, and $\boldsymbol{\gamma}$ is lower triangular.
Let $\otimes$ denote the Kronecker product. We also introduce the following notation which will be used frequently throughout the paper:
\begin{equation*}
\boldsymbol{\alpha}  \kron{\nvar} \mathbf{k} \coloneqq \bigl( \boldsymbol{\alpha} \otimes \Idvar \bigr) \, \mathbf{k}.
\end{equation*}
The Rosenbrock method \cref{eqn:ROS} is written in compact matrix notation as follows:
\begin{subequations}
\label{eqn:ROS-matrix}
\begin{align}
\label{eqn:ROS-stage-matrix}
\mathbf{k} &= h\,\fun\left( \one_s \otimes \y_{n} + \boldsymbol{\alpha} \kron{\nvar} \mathbf{k} \right)  
 + \bigl( \Idstage \otimes h\,\Jac_{n} \bigr)\, \bigl( \boldsymbol{\gamma} \kron{\nvar} \mathbf{k}  \bigr) , \\
\label{eqn:ROS-solution-matrix}
\y_{n+1} &= \y_{n} + \mathbf{b}^T \kron{\nvar} \mathbf{k},
\end{align}
where $\one_s \in \Re^s$ is a vector of ones, $\Idstage \in \Re^{s \times s}$ is the identity matrix, and
\begin{equation}
\mathbf{k} = \begin{bmatrix} k_1 \\ \vdots \\ k_s \end{bmatrix} \in \Re^{\nvar s}, \quad
\fun\left( \one_s \otimes \y_{n} + \boldsymbol{\alpha} \kron{\nvar} \mathbf{k} \right)   = \begin{bmatrix} \fun( \yy_{n}+\sum_{j}\alpha_{1,j} \,k_j ) \\ \vdots \\ \fun( \yy_{n}+\sum_{j}\alpha_{s,j} \,k_j)  \end{bmatrix} \in \Re^{\nvar s}.
\end{equation}
\end{subequations}
The Rosenbrock formula \cref{eqn:ROS-matrix} makes explicit use of the exact Jacobian, and consequently the accuracy of the method depends on the availability of the exact $\Jac_n$. In many practical cases an exact Jacobian is difficult to compute, however approximate Jacobians may be available at reasonable computational cost. Rosenbrock-W methods \citep{Steihaug_1979} maintain the accuracy of the solution when any approximation of the Jacobian is used. Specifically, an $s$-stage Rosenbrock-W method has the form \cref{eqn:ROS-matrix} but with the exact Jacobian $\Jac_n$ replaced by an arbitrary, solution-independent matrix $\Lb$  \citep[Section IV.7]{Hairer_book_II}:
\begin{subequations}
\label{eqn:ROW-matrix}
\begin{align}
\label{eqn:ROW-stage-matrix}
\mathbf{k} &= h\,\fun\left( \one_s \otimes \y_{n} + \boldsymbol{\alpha} \kron{\nvar} \mathbf{k} \right)  
 + \bigl( \Idstage \otimes h\,\Lb \bigr)\, \bigl( \boldsymbol{\gamma} \kron{\nvar} \mathbf{k}  \bigr) , \\
\label{eqn:ROW-solution-matrix}
\y_{n+1} &= \y_{n} + \mathbf{b}^T \kron{\nvar} \mathbf{k}.
\end{align}
\end{subequations}
Rosenbrock methods have received considerable attention over the years \citep{Augustin_2015_Rosenbrock,Lang_2020_Rosenbrock-overview}.
Rosen\-brock-W methods of high order have been constructed in  \citep{Rahunanthan_2010_W-methods,Rang_2005_W-methods}.
In contrast to classical interpolation/extrapolation-based multirate Rosenbrock methods~\citep{Guenther_1993_MR-ROW}, generalized multirate Rosenbrock-Wanner schemes have been introduced in~\citep{Bartel_2002_MR-W} as a special instance of partitioned Rosenbrock-W schemes. 
Matrix-free Rosen\-brock-W methods were proposed in \citep{Wensch_2005_Krylov-ROW,Schmitt_1995_W-matrix-free}, and Rosenbrock-Krylov methods that approximate the Jacobian in an Arnoldi space in \citep{Sandu_2020_BOROK,Sandu_2019_EPIRKW,Sandu_2020_ROK-adaptive,Sandu_2013_supercomputing,Sandu_2014_expK,Sandu_2014_ROK}. Application of Rosenbrock methods to parabolic partial differential equations, and the avoidance of order reduction, have been discussed in  \citep{Lubich_1995_linearlyImplicit,AlonsoMallo2002_Ros-parabolic,Schwitzer_1995_W-parabolic,Calvo_2003_order-reduction}. Linearly implicit linear multistep methods have been developed in \citep{Akrivis_2004_linearlyimplicit,Akrivis_2015_BDF-parabolic,Song_2017_linearized-BDF,Yang_2015_linearized-BDF,Yao_2017_linearized-BDF,Sandu_2020_LIMM}.

Here we consider multimethods for solving multiphysics partitioned systems \cref{eqn:additively-partitioned-ode}. Roughly speaking, multimethods allow to discretize each physical process in \cref{eqn:additively-partitioned-ode} with an appropriate numerical method; the methods are coupled appropriately such that the overall solution has the desired accuracy and stability properties. 
An example of multimethods is offered by the general-structure additive Runge--Kutta (GARK) framework, proposed in  \citep{Sandu_2015_GARK,Sandu_2016_GARK-MR}, which extends Runge--Kutta schemes to solve partitioned systems \cref{eqn:additively-partitioned-ode}. One step of a GARK method applied to the additively partitioned initial value problem \cref{eqn:additively-partitioned-ode} reads:
%
%
\begin{subequations}
\label{eqn:GARK}
\begin{align}
\label{eqn:GARK-stages}
Y^{\{q\}} &= \one_{s^{\{q\}}} \otimes \y_{n} + h \sum_{m=1}^\nparts \mathbf{A}^{\{q,m\}} \kron{\nvar} \fun^{\{m\}} (Y^{\{m\}}), \quad q=1,\ldots \nparts, \\
\label{eqn:GARK-solution}
\y_{n+1} &= \y_{n} + h\,\sum_{q=1}^\nparts \mathbf{b}^{\{q\}T} \kron{\nvar} \fun^{\{q\}} (Y^{\{q\}}).
\end{align}
\end{subequations}
Each component $\fun^{\{m\}}$ is solved with a Runge--Kutta method with $s^{\{m\}}$ stages and coefficients $(\mathbf{A}^{\{m,m\}}$, $\mathbf{b}^{\{m\}})$. The coefficients $\mathbf{A}^{\{q,m\}}$, $q \ne m$, realize the coupling among subsystems. The method \cref{eqn:GARK} builds separate stage vectors $Y^{\{m\}}$ for each component.

In this paper we construct linearly implicit multimethods that apply a possibly different Rosenbrock or Rosenbrock-W method to each component in \cref{eqn:additively-partitioned-ode}. The new family of methods, called GARK-Rosenbrock(-W), extends linearly implicit methods to solve  partitioned systems in the same way that the GARK approach \cref{eqn:GARK} extends Runge--Kutta schemes. 

The  new general GARK-Rosenbrock framework encompasses previously developed methods as particular cases, and allows one to treat them in a unified manner. The partitioned Rosenbrock-W methods, which are used in (Bartel, Günther \citep{Bartel_2002_MR-W} for deriving the order conditions of multirate W-methods, are a special case of implicit-implicit (IMIM) GARK-Rosenbrock schemes developed herein.  The Ros-ERK methods proposed by Hai and Yagi \citep{Hai_2014_ERK-ROS}, which combine explicit Runge-Kutta and Rosenbrock schemes, are a particular case of implicit-explicit (IMEX) GARK-Rosenbrock methods, and are also a particular case of methods in \citep{Bartel_2002_MR-W}. The partitioned Rosenbrock methods of Wensch et al. \citep{Wensch_1998_partitioned-ROS,Wensch_1996_partitioned-ROS} can also be cast as GARK-Rosenbrock methods (with additional special properties for solving index-3 differential algebraic equations).

The unified framework extends previous work in multiple ways. It allows to treat systems with any number of partitions, the general order conditions theory encompasses both the exact and inexact Jacobian cases, and the case where one of the partitions is very stiff (making the system an index-1 differential algebraic equation) is considered. The flexibility of the GARK-Rosenbrock framework opens the possibility to derive many new types of partitioned schemes for multiphysics systems.

The remainder of this paper is organized as follows.
\Cref{sec:Rosenbrock-partitioned} defines the new families of GARK-Rosenbrock and GARK-Rosenbrock-W methods   in the ordinary differential equation (ODE) setting.  
The order conditions theory for the new schemes is developed in \cref{sec:order-conditions} using Butcher series over special sets of trees, and linear stability is discussed in \cref{sec:linear-stability}.
 
\Cref{sec:multimethods-GARK-ROW}
constructs decoupled GARK-ROW schemes that are implicit in only one process at a time. We use  
the GARK-ROW framework to develop multimethods where each process in \cref{eqn:additively-partitioned-ode} can be solved with either an explicit Runge--Kutta, an implicit Runge--Kutta, or a Rosenbrock-W method.
Order conditions for GARK-ROS schemes applied to index-1 differential-algebraic systems are studied in \cref{sec:DAE-index-1}.
New GARK-ROW methods for practical use are proposed in \cref{sec:methods} and used for numerical experiments in \cref{sec:experiments}.
A discussion of the results in \cref{sec:discussion} concludes the paper.

\section{Partitioned Rosenbrock methods}
\label{sec:Rosenbrock-partitioned}
%
\subsection{Additively partitioned systems}

%

GARK methods \cref{eqn:GARK} extend Runge--Kutta schemes to solve partitioned systems \cref{eqn:additively-partitioned-ode}. In a similar approach, we now extend Rosenbrock methods \cref{eqn:ROS} to solve partitioned systems \cref{eqn:additively-partitioned-ode}.
Just like Rosenbrock methods are obtained by a linearization of diagonally implicit Runge--Kutta schemes, GARK-ROS methods are obtained by a linearization of  diagonally implicit GARK schemes.

\begin{definition}[GARK-ROS method]
One step of a GARK Rosenbrock   (for short, GARK-ROS)    method applied to solve the additively partitioned system \cref{eqn:additively-partitioned-ode} advances the numerical solution as follows:
\begin{subequations}
\label{eqn:GARK-ROS-scalar}
\begin{align}
\label{eqn:GARK-ROS-scalar-stages}
k_i^{\{q\}} &= h\,\fun^{\{q\}} \mleft(\y_{n} +  \sum_{m=1}^\nparts  \sum_{j=1}^{i-1} \alpha_{i,j}^{\{q,m\}}\,k_j^{\{m\}} \mright) 
 + h \, \Jac_n^{\{q\}}\,\sum_{m=1}^\nparts \sum_{j=1}^{i} \gamma_{i,j}^{\{q,m\}}\,k_j^{\{m\}} \\
 & \quad {\rm for}\quad i = 1, \ldots, s^{\{q\}},  \quad q=1, \ldots, \nparts, \nonumber\\
\label{eqn:GARK-ROS-scalar-solution}
\y_{n+1} & = \y_{n} +\sum_{q=1}^\nparts \sum_{i=1}^{s^{\{q\}}} b_{i}^{\{q\}}\,k_i^{\{q\}}.
\end{align}
\end{subequations}
The GARK-ROS scheme \cref{eqn:GARK-ROS-scalar} is written compactly in matrix notation as follows:
\begin{subequations}
\label{eqn:GARK-ROS}
\begin{align}
\label{eqn:GARK-ROS-stage}
\mathbf{k}^{\{q\}} &= h\,\fun^{\{q\}}\mleft( 
\one^{\{q\}} \otimes \y_{n}  + \sum_{m=1}^\nparts \boldsymbol{\alpha}^{\{q,m\}} \kron{\nvar}  \, \mathbf{k}^{\{m\}} \mright)    \\
\nonumber
& \quad +  \bigl( \Id_{s^{\{q\}}} \otimes h\,\Jac_n^{\{q\}} \bigr)\,\sum_{m=1}^\nparts  \boldsymbol{\gamma}^{\{q,m\}} \kron{\nvar}\, \mathbf{k}^{\{m\}}, \quad q = 1,\dots,\nparts, \\ 
\label{eqn:GARK-ROS-solution}
\y_{n+1} &= \y_{n} + \sum_{m=1}^\nparts \b^{\{m\}\,\!T} \kron{\nvar}  \mathbf{k}^{\{m\}},
\end{align}
\end{subequations} 
where we used the matrix notation \cref{eqn:ROS-matrix}. The coefficients $\boldsymbol{\alpha}^{\{q,m\}}$ are strictly lower triangular and  $\boldsymbol{\gamma}^{\{q,m\}}$ lower triangular for all $1 \le q,m \le \nparts$. 
The matrices $\Jac_n^{\{q\}} = \fun^{\{q\}}_{\y}(\y_n)$ are the Jacobians of the component functions $\fun^{\{q\}}$, evaluated at current solution $\y_n$, for each $q=1,\dots,\nparts$.
\end{definition}
The GARK-ROS scheme \cref{eqn:GARK-ROS} is characterized by the extended Butcher tableau:
\begin{equation}
\label{eqn:GARK-Rosenbrock-butcher}
\begin{butchertableau}{c|c}
\A  & \G  \\ 
\hline
\b^T &
\end{butchertableau}
~=~
\raisebox{17pt}{$
\begin{butchertableau}{ccc | ccc}
\boldsymbol{\alpha}^{\{1,1\}} &  \cdots & \boldsymbol{\alpha}^{\{1,\nparts\}} & \boldsymbol{\gamma}^{\{1,1\}} &  \ldots & \boldsymbol{\gamma}^{\{1,\nparts\}}  \\
\vdots &  \ddots & \vdots & \vdots & \ddots & \vdots \\
\boldsymbol{\alpha}^{\{\nparts,1\}} & \ \cdots & \boldsymbol{\alpha}^{\{\nparts,\nparts\}} & \boldsymbol{\gamma}^{\{\nparts,1\}} &  \ldots & \boldsymbol{\gamma}^{\{\nparts,\nparts\}} \\ 
\hline
\b^{\{1\}T} &  \cdots &\b^{\{\nparts\}T} &  &    &
\end{butchertableau}
$}.
\end{equation}

\begin{remark}[GARK-ROS scheme structure]
The GARK-ROS scheme \linebreak \cref{eqn:GARK-ROS} has the following characteristics:
\begin{itemize}
\item A different increment vector $\mathbf{k}^{\{q\}} \in \Re ^{\nvar s}$ is constructed for each component $q=1,\dots,\nparts$.
\item Computation of the increment $\mathbf{k}^{\{q\}}$ uses only evaluations of the corresponding component function $\fun^{\{q\}}$. The argument at which $\fun^{\{q\}}$ is   evaluated   is  constructed using a linear combination of all increments $\mathbf{k}^{\{m\}}$ for $m=1,\dots,\nparts$.
\item Computation of the increment $\mathbf{k}^{\{q\}}$ involves linear combinations of increments $\mathbf{k}^{\{m\}}$ for $m=1,\dots,\nparts$, multiplied by the Jacobian $\Jac^{\{q\}}$ of the corresponding component function. Therefore the calculation of increments involves the solution of linear systems.
\item For all $\gamma_{i,j}^{\{q,m\}}=0$, the scheme \cref{eqn:GARK-ROS} reduces to an explicit GARK method.
\item If $\gamma_{i,j}^{\{q,m\}}=0$ for all $m>q$ holds, all increments can be computed recursively: $k_1^{\{1\}}, \ldots, k_1^{\{\nparts\}}, k_2^{\{1\}},\ldots, k_{s^{\{\nparts\}}}^{\{\nparts\}}$.
\end{itemize}
\end{remark}

\begin{definition}[GARK-ROW method]
One step of a GARK Rosenbrock-W (for short, GARK-ROW) method applied to solve the additively partitioned system \cref{eqn:additively-partitioned-ode} advances the numerical solution as follows:
\begin{subequations}
\label{eqn:GARK-ROW}
\begin{align}
\label{eqn:GARK-ROW-stage}
\mathbf{k}^{\{q\}} &= h\,\fun^{\{q\}}\mleft( 
\one^{\{q\}} \otimes \y_{n}  + \sum_{m=1}^\nparts \boldsymbol{\alpha}^{\{q,m\}} \kron{\nvar}  \, \mathbf{k}^{\{m\}} \mright)    \\
\nonumber
& \quad +  \bigl( \Id_{s^{\{q\}}} \otimes h\,\Lb^{\{q\}} \bigr)\,\sum_{m=1}^\nparts  \boldsymbol{\gamma}^{\{q,m\}} \kron{\nvar}\, \mathbf{k}^{\{m\}}, \quad q = 1,\dots,\nparts, \\ 
\label{eqn:GARK-ROW-solution}
\y_{n+1} &= \y_{n} + \sum_{m=1}^\nparts \b^{\{m\}\,\!T} \kron{\nvar}  \mathbf{k}^{\{m\}}.
\end{align}
\end{subequations} 
where $\Lb^{\{q\}}$ are arbitrary matrix approximations to component function Jacobians $\fun^{\{q\}}_{\y}(\y_n)$, for each $q=1,\dots,\nparts$.
\end{definition}

\subsection{Component partitioned systems}
%
Consider the partitioned system:
\begin{equation}
\label{eqn:component-partitioned-ode}
\frac{d\yy^{\{q\}}}{dt} = \fun^{\{q\}}\mleft( \mathbf{y}^{\{1\}}, \cdots, \mathbf{y}^{\{\nparts\}} \mright), ~~
\mathbf{y}^{\{q\}} \in \Re^{\nvar^{\{q\}}}, ~~ q=1,\dots,\nparts 
\end{equation}
with $\sum_{q=1}^{\nparts} \nvar^{\{q\}} = \nvar.$
The Jacobian of each component function $\fun^{\{q\}}$ with respect to each component vector is approximated by:
\begin{equation*}
\frac{ \partial\,\fun^{\{q\}}}{\partial\,\mathbf{y}^{\{m\}}} = (\fun_\y)^{\{q,m\}} \approx \Lb^{\{q,m\}} \in \Re^{\nvar^{\{q\}} \times \nvar^{\{m\}}}.
\end{equation*}
The GARK-ROW scheme \cref{eqn:GARK-ROW} applied to a component split system \cref{eqn:component-partitioned-ode} reads:
\ifreport
\begin{subequations}
\label{eqn:GARK-ROW-scalar-component}
\begin{align}
\label{eqn:GARK-ROW-scalar-stages-component}
k_i^{\{q\}} &= h\,\fun^{\{q\}} \mleft(\y_n^{\{1\}} + \sum_{j=1}^{i-1 } \alpha_{i,j}^{\{q,1\}}\,k_j^{\{1\}},\cdots, 
y_n^{\{\nparts\}} + \sum_{j=1}^{i-1 } \alpha_{i,j}^{\{q,\nparts\}}\,k_j^{\{\nparts\}}\mright) \\
\nonumber
&\quad + h \, \sum_{m=1}^\nparts \sum_{j=1}^i \gamma_{i,j}^{\{q,m\}}\,\Lb^{\{q,m\}}\, k_j^{\{m\}},\quad i=1,\dots,s, \\
\label{eqn:GARK-ROW-scalar-solution-component}
y_{n+1}^{\{q\}} & = y_n^{\{q\}} + \sum_{i=1}^{s^{\{q\}}} b_{i}^{\{q\}}\,k_i^{\{q\}}, \quad q = 1,\dots,\nparts.
\end{align}
\end{subequations}
In matrix notation the GARK-ROW scheme \cref{eqn:GARK-ROW-scalar-component} applied to the component split system \cref{eqn:component-partitioned-ode} reads:
\fi
\begin{subequations}
\label{eqn:GARK-ROW-component}
\begin{align}
\mathbf{Y}^{\{q,m\}} &= \one_{s^{\{q\}}} \otimes  \y^{\{m\}}_{n}  + \bigl( \boldsymbol{\alpha}^{\{q,m\}} \otimes \mathbf{I}_{\nvar^{\{m\}}} \bigr) \, \mathbf{k}^{\{m\}} \in \Re^{\nvar^{\{m\}} s^{\{q\}}}, \\
\label{eqn:GARK-ROW-stage-component}
\mathbf{k}^{\{q\}} &= h\,\fun^{\{q\}}\mleft( \mathbf{Y}^{\{q,1\}},\cdots, \mathbf{Y}^{\{q,\nparts\}} \mright)    + h\, \sum_{m=1}^\nparts  \Bigl( \boldsymbol{\gamma}^{\{q,m\}} \otimes \mathbf{L}^{\{q,m\}} \Bigr) \, \mathbf{k}^{\{m\}},  \\ 
\label{eqn:GARK-ROW-solution-component}
\y^{\{q\}}_{n+1} &= \y^{\{q\}}_{n} + \bigl( \b^{\{q\}\,\!T} \otimes \mathbf{I}_{\nvar^{\{q\}}} \bigr) \, \mathbf{k}^{\{q\}}, \quad q = 1,\dots,\nparts.
\end{align}
\end{subequations} 

\begin{remark}
The GARK-ROS scheme \cref{eqn:GARK-ROS} applied to a component split system \cref{eqn:component-partitioned-ode} has the form \cref{eqn:GARK-ROW-component}, where each matrix equals the corresponding sub-Jacobian $\Lb^{\{q,m\}} = \partial\,\fun^{\{q\}}/\partial\,\mathbf{y}^{\{m\}}(\y_n)$. Thus component partitioned systems are a special case of additively partitioned systems. 
\end{remark}

\section{Order conditions}
\label{sec:order-conditions}
%
We develop the order conditions theory for additively partitioned systems \cref{eqn:additively-partitioned-ode}. These order conditions remain valid for component partitioned systems \cref{eqn:component-partitioned-ode} as well.

\subsection{Multicolored trees and NB-series}

We  recall the set of $\NT$ trees 
\citep{SanzSerna_1997_symplectic} which provide a generalization of Butcher trees for partitioned systems.
\begin{definition}
The set $\NT$ consists of rooted trees with round ($\mround{m}$) vertices, each colored in one of the distinct $m=1,\dots,\nparts$ colors. Here nodes of color $m$ correspond to derivatives of the component function $\fun^{\{m\}}$ of the partitioned system \cref{eqn:additively-partitioned-ode}. 
\end{definition}
We now introduce the set of trees that represent the GARK-ROW numerical solution.
\begin{definition}
The set $\NTW$ consists of rooted trees with both square ($\msquare{m}$) and round ($\mround{m}$) vertices, each colored in one of the distinct $m=1,\dots,\nparts$ colors. Square nodes have a single child, and there are no square leaves. Each color corresponds to a different component of the partitioned system. For our purpose, round nodes ($\mround{m}$) represent derivatives of the component function $\fun^{\{m\}}$, and square nodes ($\msquare{m}$) to the action of the partition's approximate Jacobian matrix $\mathbf{L}^{\{m\}}$. 
\end{definition}
\begin{remark}
\label{rem:nt-in-ntw}
Clearly $\NT \subset \NTW$. The following properties discussed for $\NTW$ are applicable to $\NT$ as well.
\end{remark}
The empty $\NTW$ tree is denoted by $\emptyset$. The $\NTW$ tree with a single vertex of color $m$ is denoted by $\tau_{\mround{m}}$. We denote by $\tree = [ \tree_1 \dots \tree_L ]_{\mround{m}} \in$ $\NTW$ the new tree obtained by joining $\tree_1, \dots, \tree_L \in$ $\NTW$ with a root of color $m$ (i.e., attaching each of the trees directly to the root, which will have $L$ children). 
We denote by $\tree = [ \tree_1  ]_{\msquare{m}} \in$ $\NTW$ the new tree obtained by appending to $\tree_1 \in$ $\NT$ a square root of color $m$.

Similar to regular Butcher trees, the order $\rho(\tree)$ is the number of nodes of $\tree \in \NTW$.
The density $\gamma(\tree)$ and the number of symmetries $\sigma(\tree)$ are defined recursively by
\begin{alignat*}{5}
\gamma(\emptyset) & = 1; \quad & \gamma\left(\tau_{\mround{m}}\right) &= 1; \quad &
\gamma(\tree) &= 
\begin{cases} 
\rho(\tree)\, \gamma(\tree_1) \cdots \gamma(\tree_L), & \textnormal{for}~\tree=[\tree_1,\ldots,\tree_L]_{\mround{m}}, \\
\rho(\tree)\, \gamma(\tree_1), & \textnormal{for}~\tree=[\tree_1]_{\msquare{m}}, \\
\end{cases}
\\
\sigma(\emptyset) &= 1; \quad & \sigma\left(\tau_{\mround{m}}\right) &= 1; \quad &
\sigma(\tree) &= 
\begin{cases} 
\prod_{i=1}^L m_i!\, \sigma(\tree_i)^{m_i}, & 
\textnormal{for}~\tree=[\tree_1^{m_1},\ldots,\tree_L^{m_L}]_{\mround{m}}, \\
\sigma(\tree_1), &  
\textnormal{for}~\tree=[\tree_1]_{\msquare{m}},
\end{cases}
\end{alignat*}
with $\tree_l^{m_l}$ meaning that the tree $\tree_l$ has been attached $m_l$ times to the root $\mround{m}$.

\begin{definition}[Elementary differentials over $\NTW$]
\label{def:elementary-differentials}
An elementary differential $F(\tree)(\cdot) : \Re^\nvar \to \Re^\nvar$ is associated to each tree $\tree \in \NTW$.
Using tensorial notation, the elementary differentials are defined recursively as follows:
\begin{equation}
\label{eqn:elementary-differentials}
F(\tree)(\y_\ast) = \begin{cases}
0, & \textnormal{for}~\tree = \emptyset; \\
\fun^{\{m\}}(\y_\ast), & \textnormal{for}~\tree = \tau_{\mround{m}}; \\
\frac{d^L \fun^{\{m\}}}{d\y^L}(\y_\ast)\Bigl( F(\tree_1)(\y_\ast), \dots, F(\tree_L)(\y_\ast)\Bigr), & \textnormal{for}~\tree = [ \tree_1 \dots \tree_L ]_{\mround{m}}; \\
\Lb^{\{m\}} \cdot F(\tree_1)(\y_\ast) & \textnormal{for}~\tree = [ \tree_1 ]_{\msquare{m}}, ~   \rho(\tree_1)   
\ge 1.
\end{cases}
\end{equation}
The second argument of the elementary differential is a vector $\y_* \in \Re^\nvar$ which represents the argument at which all the function derivatives are evaluated.
\end{definition}

We extend the Butcher series (B-series) to the sets $\NT$ and $\NTW$.
\begin{definition}
An NB-series is a formal expansion in powers of the step size $h$
\begin{equation}
\label{eqn:NB-series-definition} 
\textnormal{NB}(\mathfrak{c},\y_\ast) \coloneqq \sum_{\tree \in \NTW} \frac{h^{\rho(\tree)}}{\sigma(\tree)}\, \mathfrak{c}(\tree)\, F(\tree)(\y_\ast)\,,
\end{equation}
where the summation is carried out over elements of a set of rooted trees. Each term consists of a weighted elementary differential \cref{eqn:elementary-differentials}. Here we consider summation over $\NTW$, with $\mathfrak{c} : \NTW \to \Re$ a mapping that assigns a real number  to each tree.  Per \cref{rem:nt-in-ntw} an NB-series over $\NT$ has the form \cref{eqn:NB-series-definition} with $\mathfrak{c}(\tree)=0$ for any $\tree \in \NTW \backslash \NT$.
\end{definition}

\begin{lemma}
The exact solution of \cref{eqn:additively-partitioned-ode} is represented by the NB-series \citep{SanzSerna_1997_symplectic}
\begin{equation}
\label{eqn:NB-exact}
\y(t+h)=\textnormal{NB}(\mathfrak{c},\y(t)) \quad \textnormal{with} \quad \mathfrak{c}(\tree)= \begin{cases}
\frac{1}{\gamma(\tree)}, & \textnormal{for}~\tree \in \NT, \\
0, & \textnormal{for}~\tree \in \NTW \backslash \NT.
\end{cases}
\end{equation}
\end{lemma}

We next provide several results that will prove useful to derive the order conditions of partitioned Rosenbrock methods.

\begin{theorem}[Function of NB-series \citep{Sandu_2015_GARK}]
\label{thm:F-of-NBseries}
A component function applied to an NB-series \cref{eqn:NB-series-definition} with $\mathfrak{a}(\emptyset)=1$ is also an NB-series, 
\begin{equation*}
\begin{split}
h\, \fun^{\{m\}}\big( \textnormal{NB}(\mathfrak{a},\y_{n}) \big) 
=\textnormal{NB}\bigl( (\mathrm{D}^{\{m\}} \mathfrak{a}), \y_{n}\bigr),
\end{split}
\end{equation*}
characterized by the coefficients:
\begin{equation}
\label{eqn:B-series-function}
(\mathrm{D}^{\{m\}} \mathfrak{a})(\tree) = \begin{cases}
0, & \textnormal{for}~\tree=\emptyset,\\
1, & \textnormal{for}~\tree  = \tau_{\mround{m}},  \\
\prod_{\ell=1}^{L} \mathfrak{a}(\tree_\ell) &
\textnormal{for}~\tree=[\tree_1,\ldots,\tree_L]_{\mround{m}},~L\ge 1, \\
0, & \textnormal{otherwise}.
\end{cases}
\end{equation}
\end{theorem}

\begin{theorem}[Jacobian times NB-series]
\label{thm:J-times-NBseries}
A Jacobian matrix times an NB-series \cref{eqn:NB-series-definition} with $\mathfrak{a}(\emptyset)=0$ is also an NB-series, 
\begin{equation*}
\begin{split}
h\, \Jac^{\{m\}}_{n} \cdot \big( \textnormal{NB}(\mathfrak{a},\y_{n}) \big) 
=\textnormal{NB}\bigl( (\mathrm{J}^{\{m\}} \mathfrak{a}), \y_{n}\bigr),
\end{split}
\end{equation*}
characterized by the coefficients:
\begin{equation}
\label{eqn:B-series-jacobian}
(\mathrm{J}^{\{m\}} \mathfrak{a})(\tree) = \begin{cases}
 \mathfrak{a}(\utree), & \textnormal{for}~\tree=[\utree]_{\mround{m}},\\
0, & \textnormal{otherwise}.
\end{cases}
\end{equation}
\end{theorem}

\begin{proof}
We consider the Jacobian matrix times the series:
\begin{equation*}
\begin{split}
h\, \Jac^{\{m\}}_{n} \cdot \big( \textnormal{NB}(\mathfrak{a},\y_{n}) \big) 
&=\sum_{\tree \in \NTW} \mathfrak{a}(\tree)\, \frac{h^{\rho(\tree)+1}}{\sigma(\tree)}\,\fun^{\{m\}}_\y(\y_{n}) \, F(\tree)\bigl(\y_{n}\bigr). 
\end{split}
\end{equation*}
This expression involves elementary differentials $\fun_\y \cdot F(\tree)$, and we note that:
\begin{equation*}
\fun^{\{m\}}_\y(\y_{n}) \cdot F(\tree)\bigl(\y_{n}\bigr) = F([\tree]_{\mround{m}})\bigl(\y_{n}\bigr),
\end{equation*}
and that $\rho([\tree]_{\mround{m}}) = \rho(\tree)+1$ and $\sigma([\tree]_{\mround{m}}) = \sigma(\tree)$, which leads to \cref{eqn:B-series-jacobian}.
\end{proof}

\begin{theorem}[Jacobian approximation times NB-series]
\label{thm:L-times-NBseries}
A Jacobian approximation matrix times an NB-series \cref{eqn:NB-series-definition} with $\mathfrak{a}(\emptyset)=0$ is also an NB-series, 
\begin{equation*}
h\, \Lb^{\{m\}}_{n} \cdot \big( \textnormal{NB}(\mathfrak{a},\y_{n}) \big) 
=\textnormal{NB}\bigl( (\mathrm{L}^{\{m\}} \mathfrak{a}), \y_{n}\bigr),
\end{equation*}
characterized by the coefficients:
\begin{equation}
\label{eqn:B-series-L}
(\mathrm{L}^{\{m\}} \mathfrak{a})(\tree) = \begin{cases}
 \mathfrak{a}(\utree), & \textnormal{for}~\tree=[\utree]_{\msquare{m}},\\
0, & \textnormal{otherwise}.
\end{cases}
\end{equation}
\end{theorem}

\begin{proof}
Similar to the proof of \cref{thm:J-times-NBseries}.
\end{proof}

\subsection{GARK-ROS order conditions}
We represent the stage vectors and numerical solutions of GARK-ROS methods \cref{eqn:GARK-ROS} as NB-series \cref{eqn:NB-series-definition} over $\NTW$: 
\begin{equation}
\label{eqn:GARK-ROS-Bseries} 
\mathbf{k}^{\{q\}} = \textnormal{NB}\left( \boldsymbol{\theta}^{\{q\}}, \y_{n} \right) \in \Re^{s^{\{q\}}}, \quad 
\y_{n+1} = \textnormal{NB}\left( \boldsymbol{\phi}, \y_{n} \right) \in \Re.
\end{equation}
Insert \cref{eqn:GARK-ROS-Bseries} into the stage equations \cref{eqn:GARK-ROS-stage}
\ifreport
\begin{equation*}
\begin{split}
\textnormal{NB}\left( \boldsymbol{\theta}^{\{q\}}, \y_{n} \right) &= h\,\fun^{\{q\}}\mleft( 
\one_s + \sum_{m=1}^\nparts \boldsymbol{\alpha}^{\{q,m\}} \, \textnormal{NB}\left( \boldsymbol{\theta}^{\{m\}}, \y_{n} \right) \mright) \\
& \quad + h\,\Jac_n^{\{q\}} \,\sum_{m=1}^\nparts  \boldsymbol{\gamma}^{\{q,m\}} \, \textnormal{NB}\left( \boldsymbol{\theta}^{\{m\}}, \y_{n} \right),
\end{split}
\end{equation*}
\fi
and apply \cref{thm:F-of-NBseries} and \cref{thm:J-times-NBseries} to obtain:
\begin{equation*}
\begin{split}
\boldsymbol{\theta}^{\{q\}}(\tree) &= {\left(\mathrm{D}^{\{q\}} \sum_{m=1}^\nparts \boldsymbol{\alpha}^{\{q,m\}}\,\boldsymbol{\theta}^{\{m\}}\right)}(\tree) + \sum_{m=1}^\nparts \boldsymbol{\gamma}^{\{q,m\}}\,\left(\mathrm{J}^{\{q\}}\boldsymbol{\theta}^{\{m\}}\right)(\tree).
\end{split}
\end{equation*}
This leads to the following recurrence on stage vectors NB-series coefficients \cref{eqn:GARK-ROS-Bseries}:
\begin{equation}
\label{eqn:GARK-ROS-K-series}
\boldsymbol{\theta}^{\{q\}}(\tree) = \begin{cases}
0, & \tree = \emptyset, \\
1, & \tree  = \tau_{\mround{q}},  \\
\Times_{\ell=1}^{L} \left(\sum_{m=1}^\nparts \boldsymbol{\alpha}^{\{q,m\}}\,\boldsymbol{\theta}^{\{m\}}(\tree_\ell)\right), &
\textnormal{for}~\tree=[\tree_1,\ldots,\tree_L]_{\mround{q}},~L\ge 2, \\
\sum_{m=1}^\nparts \boldsymbol{\beta}^{\{q,m\}}\,\boldsymbol{\theta}^{\{m\}}(\tree_{1}), & \textnormal{for}~\tree=[\tree_1]_{\mround{q}}, \\
0, & \textnormal{when root}(\tree) \ne \mround{q}.
\end{cases}
\end{equation}
We denote by $\Times$ the {\it element-by-element} product of $s$-dimensional vectors.
Note that in sums of the form $\sum_{m=1}^\nparts \boldsymbol{\alpha}^{\{q,m\}}\,\boldsymbol{\theta}^{\{m\}}(\tree)$ and $\sum_{m=1}^\nparts \boldsymbol{\beta}^{\{q,m\}}\,\boldsymbol{\theta}^{\{m\}}(\tree)$ at most a single term is nonzero, namely, the one with $m = n$ when  root$(\tree) = \mround{n}$. The recurrence \cref{eqn:GARK-ROS-K-series} only builds terms corresponding to trees in $\NT$; consequently, $\boldsymbol{\theta}^{\{q\}}(\tree)=0$ for $\tree \in \NTW \backslash \NT$.

Inserting \cref{eqn:GARK-ROS-Bseries} into the solution equations \cref{eqn:GARK-ROS-solution} leads to the following B-series coefficients of the numerical solution:
\begin{equation}
\label{eqn:GARK-ROS-Y-series}
\begin{split}
\ifreport
\textnormal{NB}\left( \boldsymbol{\phi}, \y_{n} \right) &= 1 + \sum_{m=1}^\nparts \b^{\{m\}\,\!T} \textnormal{NB}\left( \boldsymbol{\theta}^{\{m\}}, \y_{n} \right) \\
\Rightarrow \quad 
\fi
\boldsymbol{\phi}(\tree) &= \begin{cases} 1, & \tree = \emptyset, \\
 \sum_{m=1}^\nparts \b^{\{m\}\,\!T} \boldsymbol{\theta}^{\{m\}}(\tree), & \tree \in \NT  \backslash \{ \emptyset \}, \\
 0, & \tree \in \NTW \backslash \NT.
 \end{cases}
\end{split}
\end{equation}
A comparison of the numerical solution \cref{eqn:GARK-ROS-Y-series} with the exact solution \cref{eqn:NB-exact} leads to the following result.
\begin{theorem}[GARK-ROS order conditions]
The GARK-ROS method \cref{eqn:GARK-ROS} has order of consistency $p$ iff
\begin{equation*}
\sum_{m=1}^\nparts \b^{\{m\}\,\!T} \boldsymbol{\theta}^{\{m\}}(\tree) = \frac{1}{\gamma(\tree)} \quad 
\textnormal{for } \tree \in \NT \textnormal{  with } 1 \le \rho(\tree) \le p.
\end{equation*}
\end{theorem}

The procedure to generate the order conditions for GARK-ROS methods using the recurrence \cref{eqn:GARK-ROS-K-series} is illustrated in  \cref{tab:GARK_ROS-trees_1-3_numeric}. The process is as follows:
\begin{itemize}
\item The root of color $m$ is labelled $\mathbf{b}^{\{m\}T}$.
\item A single sibling of color $m$ (its parent of color $q$ has one child) is labelled $\boldsymbol{\beta}^{\{q,m\}}$.
\item A node of color $m$ with multiple siblings (its parent of color $q$ has  multiple children) is labelled $\boldsymbol{\alpha}^{\{q,m\}}$. 
\item The result of each subtree is an $s$-dimensional vector of NB-series coefficients.
\item The leaves build their vector by multiplying their label by a vector of ones. 
\item A node (except the leaves) takes the element-wise product of the vectors of its children, then multiplies the result by its label. 
\end{itemize} 
 
\begin{table}
\centering
\begin{tabular}{|| C{2em} || C{7.5em} || C{7em}|| C{10.5em}  ||  C{2em} ||}
\hline\hline
$\tree$ &  Labels & $F(\tree)$ & $\phi(\tree)$ & $\gamma(\tree)$ \\
\hline\hline
$\tree_1$ &
\begin{tikzpicture}[grow=up, level distance=0.8cm, scale=0.75]\Tree [.\node[circle,draw,thick,opacity=1,fill=none,inner sep=1,minimum size=4,label=right:$\mathbf{b}^{\{m\}T}$](b){m};  ]\end{tikzpicture} 
&  $\fun^{\{m\}}$ & $\mathbf{b}^{\{m\}T} \, \one^{\{m\}}$ & 1 \\ 
\hline\hline
$\tree_2$ &
\begin{tikzpicture}[grow=up, level distance=0.8cm, scale=0.75]\Tree [.\node[circle,draw,thick,opacity=1,fill=none,inner sep=1,minimum size=4,label=right:$\mathbf{b}^{\{m\}T}$](j){m}; [.\node[circle,draw,thick,opacity=1,fill=none,inner sep=1,minimum size=4,label=right:$\boldsymbol{\beta}^{\{m,n\}}$](k){n};  ] ]\end{tikzpicture} 
& $\fun_\y^{\{m\}}\,\fun^{\{n\}}$ &   $\mathbf{b}^{\{m\}T} \, \boldsymbol{\beta}^{\{m,n\}} \, \one^{\{n\}}$ & 2 \\ 
\hline\hline
$\tree_{3,1}$ &
\begin{tikzpicture}[grow=up, level distance=0.8cm, scale=0.75]\Tree [.\node[circle,draw,thick,opacity=1,fill=none,inner sep=1,minimum size=4,label=right:$\mathbf{b}^{\{m\}T}$](j){m}; [.\node[circle,draw,thick,opacity=1,fill=none,inner sep=1,minimum size=4,label=right:$\boldsymbol{\alpha}^{\{m,p\}}$](l){p};  ][.\node[circle,draw,thick,opacity=1,fill=none,inner sep=1,minimum size=4,label=right:$\boldsymbol{\alpha}^{\{m,n\}}$](k){n};  ] ]\end{tikzpicture} &   
$\fun_{\y,\y}^{\{m\}}\,(\fun^{\{n\}},\fun^{\{p\}}) $ & $\begin{array}{c} \mathbf{b}^{\{m\}T} \, \bigl( (\boldsymbol{\alpha}^{\{m,n\}}\, \one^{\{n\}}) \\ \times (\boldsymbol{\alpha}^{\{m,p\}}\, \one^{\{p\}}) \bigr) \end{array} $ & 3 \\ 
\hline
$\tree_{3,2}$ &
\begin{tikzpicture}[grow=up, level distance=0.8cm, scale=0.75]\Tree [.\node[circle,draw,thick,opacity=1,fill=none,inner sep=1,minimum size=4,label=right:$\mathbf{b}^{\{m\}T}$](j){m}; [.\node[circle,draw,thick,opacity=1,fill=none,inner sep=1,minimum size=4,label=right:$\boldsymbol{\beta}^{\{m,n\}}$](k){n}; [.\node[circle,draw,thick,opacity=1,fill=none,inner sep=1,minimum size=4,label=right:$\boldsymbol{\beta}^{\{n,p\}}$](l){p};  ] ] ]\end{tikzpicture} &
$\fun_\y^{\{m\}}\,\fun_\y^{\{n\}}\,\fun^{\{p\}}$ & $\begin{array}{c} \mathbf{b}^{\{m\}T}  \, \boldsymbol{\beta}^{\{m,n\}} \\ \cdot \boldsymbol{\beta}^{\{n,p\}}\,\one^{\{p\}} \end{array}$ & 6 \\ 
\hline\hline
\end{tabular}
\caption{$\NT$ trees of orders 1 to 3 for the GARK-ROS numerical solution. The root of color $m$ is labelled $\mathbf{b}^{\{m\}T}$.
Single siblings are labelled $\boldsymbol{\beta}$, vertices that have multiple siblings are labelled $\boldsymbol{\alpha}$, and each node label is superscripted by a pair of indices $\{q,m\}$, where $m$ is the color of the node  and  $q$ the color of its parent.}
\label{tab:GARK_ROS-trees_1-3_numeric}
\end{table}


We note that each node (except the roots) carries a label with two indices, first the color of its parent, followed by its own color. Moreover, if all the nodes have the same color then $\NT$ is the set of T-trees, and the GARK-ROS order conditions give the Rosenbrock order conditions. These observations lead to the following result.
\begin{theorem}[GARK-ROS order conditions]
The GARK-ROS order conditions \cref{eqn:GARK-ROS}  are the same as the Rosenbrock order conditions \cref{eqn:ROS}, except that the method coefficients are labelled according to node colors. In the order conditions, in each sequence of matrix multiplies, the color indices are compatible according to matrix multiplication rules.
\end{theorem}

Let $\one^{\{n\}} \in \Re^{s^{\{n\}}}$ be a vector of ones.   For brevity we also define the  vectors:
\begin{equation}
\label{eqn:terminal-vectors}
\begin{split}
\mathbf{c}^{\{m,n\}} &\coloneqq \boldsymbol{\alpha}^{\{m,n\}} \,  \one^{\{n\}}, \quad
\mathbf{g}^{\{m,n\}} \coloneqq \boldsymbol{\gamma}^{\{m,n\}} \,  \one^{\{n\}}, \\
\mathbf{e}^{\{m,n\}} &\coloneqq \boldsymbol{\beta}^{\{m,n\}} \, \one^{\{n\}}
= \mathbf{c}^{\{m,n\}} + \mathbf{g}^{\{m,n\}}.
\end{split}
\end{equation}
The GARK-ROS order four conditions read: 
\begin{subequations}
\label{eqn:GARK-ROS-order-conditions}
\begin{align}
\label{eqn:GARK-ROS-order1-conditions}
\textnormal{order 1:   } &\begin{cases}
\mathbf{b}^{\{m\}\,T} \, \one^{\{m\}} = 1, \end{cases} \textnormal{for } m = 1,\dots,\nparts; \\
\label{eqn:GARK-ROS-order2-conditions}
\textnormal{order 2:   } &
\begin{cases}
\mathbf{b}^{\{m\}\,T} \, \mathbf{e}^{\{m,n\}}  = \frac{1}{2}, \end{cases}  \textnormal{for } m,n = 1,\dots,\nparts; \\
\label{eqn:GARK-ROS-order3-conditions}
\textnormal{order 3:   } &
\begin{cases}
\mathbf{b}^{\{m\}\,T} \, \bigl( \mathbf{c}^{\{m,n\}} \times \mathbf{c}^{\{m,p\}} \bigr) = \frac{1}{3}, \\
\mathbf{b}^{\{m\}\,T} \, \boldsymbol{\beta}^{\{m,n\}}\, \mathbf{e}^{\{n,p\}} = \frac{1}{6},
\end{cases}  \textnormal{for } m,n,p = 1,\dots,\nparts; \\
\label{eqn:GARK-ROS-order4-conditions}
\textnormal{order 4:   } &
\begin{cases}
\mathbf{b}^{\{m\}\,T} \,\left( \mathbf{c}^{\{m,n\}} \times  \mathbf{c}^{\{m,p\}}  \times  \mathbf{c}^{\{m,q\}}  \right) = \frac{1}{4}, \\
\mathbf{b}^{\{m\}\,T} \, \big((\boldsymbol{\alpha}^{\{m,n\}}\,\mathbf{e}^{\{n,p\}}) \times \mathbf{c}^{\{m,p\}} \big) = \frac{1}{8}, \\
\mathbf{b}^{\{m\}\,T} \, \boldsymbol{\beta}^{\{m,n\}}\, \left( \mathbf{c}^{\{n,p\}} \times \mathbf{c}^{\{n,q\}}  \right) = \frac{1}{12}, \\ 
\mathbf{b}^{\{m\}\,T} \,\boldsymbol{\beta}^{\{m,n\}}\,\boldsymbol{\beta}^{\{n,p\}}\,\mathbf{e}^{\{p,q\}} =  \frac{1}{24}, \\
 \textnormal{for } m,n,p,q = 1,\dots,\nparts.
\end{cases} 
\end{align}
\end{subequations}
%

\subsection{GARK-ROW order conditions}

We represent the stage vectors and numerical solutions of GARK-ROW    methods    \cref{eqn:GARK-ROW}   
as NB-series \cref{eqn:NB-series-definition} over $\NTW$: 
\begin{equation}
\label{eqn:GARK-ROW-Bseries} 
\mathbf{k}^{\{q\}} = \textnormal{NB}\left( \boldsymbol{\theta}^{\{q\}}, \y_{n} \right) \in \Re^s, \quad 
\y_{n+1} = \textnormal{NB}\left( \boldsymbol{\phi}, \y_{n} \right) \in \Re.
\end{equation}
Insert \cref{eqn:GARK-ROW-Bseries} into the stage equations \cref{eqn:GARK-ROW-stage}, and apply \cref{thm:F-of-NBseries} and \cref{thm:L-times-NBseries} to obtain:
\begin{equation*}
\begin{split}
\boldsymbol{\theta}^{\{q\}}(\tree) &= {\left(\mathrm{D}^{\{q\}} \sum_{m=1}^\nparts \boldsymbol{\alpha}^{\{q,m\}}\,\boldsymbol{\theta}^{\{m\}}\right)}(\tree) + \sum_{m=1}^\nparts \boldsymbol{\gamma}^{\{q,m\}}\,\left(\mathrm{L}^{\{q\}}\boldsymbol{\theta}^{\{m\}}\right)(\tree).
\end{split}
\end{equation*}
This leads to the following recurrence on NB-series coefficients:
\begin{equation*}
\begin{split}
\boldsymbol{\theta}^{\{q\}}(\tree) &= \begin{cases}
0, & \tree = \emptyset, \\
1, & \tree  = \tau_{\mround{q}},  \\
\Times_{\ell=1}^{L} \left(\sum_{m=1}^\nparts \boldsymbol{\alpha}^{\{q,m\}}\,\boldsymbol{\theta}^{\{m\}}(\tree_\ell)\right), &
	\textnormal{for}~\tree=[\tree_1,\ldots,\tree_L]_{\mround{q}},~L\ge 1, \\
\sum_{m=1}^\nparts \boldsymbol{\gamma}^{\{q,m\}}\,\boldsymbol{\theta}^{\{m\}}(\tree_{1}), & 
	\textnormal{for}~\tree=[\tree_1]_{\msquare{q}}, \\
0, &  \textnormal{when root}(\tree) \not\in \{ \mround{q}, \msquare{q} \}.
\end{cases}
\end{split}
\end{equation*}
Note that in sums of the form $\sum_{m=1}^\nparts \boldsymbol{\gamma}^{\{q,m\}}\,\boldsymbol{\theta}^{\{m\}}(\tree_{1})$ a single term is nonzero, namely, the one with $m$ equal the color of the root of $\tree_1$.

Inserting \cref{eqn:GARK-ROW-Bseries} into the solution equations \cref{eqn:GARK-ROW-solution} leads to an NB-series representation of the numerical solution given by \cref{eqn:GARK-ROS-Y-series}. Equating the terms of the numerical solution NB-series with those of the exact solution \cref{eqn:NB-exact} leads to the following order conditions theorem.
\begin{theorem}[GARK-ROW order conditions]
The GARK-ROW method \cref{eqn:GARK-ROW} has order $p$ iff:
\begin{equation*}
\boldsymbol{\phi}(\tree) = \begin{cases}
\frac{1}{\gamma(\tree)}, & \textnormal{for}~ \tree \in \NT,\\
0, & \textnormal{for}~ \tree \in \NTW \backslash \NT,
\end{cases}
\quad \textnormal{for } \tree \in \NTW \textnormal{  with } 1 \le \rho(\tree) \le p.
\end{equation*}
\end{theorem}

The procedure to generate the order conditions for GARK-ROS methods using the recurrence \cref{eqn:GARK-ROS-K-series} is illustrated in  \cref{tab:GARK-ROW-trees_1-3_numeric}. The process is as follows:
\begin{itemize}
\item Roots of color $q$ are labelled $\mathbf{b}^{\{q\}}$;
\item Nodes of color $m$ with a round parent of color $q$ are labelled $\boldsymbol{\alpha}^{\{q,m\}}$;
\item Nodes of color $m$ with a square parent of color $q$ are labelled $\boldsymbol{\gamma}^{\{q,m\}}$;
\item The result of each subtree is an $s$-dimensional vector of NB-series coefficients. Obtaining these coefficients is done starting from the leaves and working toward the root, as discussed for GARK-ROS methods.
\end{itemize}

We note that each node (except the roots) carries a label with two indices, first the color of its parent, followed by its own color. Moreover, if all the nodes have the same color then $\NTW$ is the set of TW-trees, and the GARK-ROW order conditions give the Rosenbrock-W order conditions. We have the following result.
\begin{theorem}[GARK-ROW order conditions]
The GARK-ROW order conditions \cref{eqn:GARK-ROW}  are the same as the Rosenbrock-W order conditions \cref{eqn:ROW-matrix}, except that the method coefficients are labelled according to node colors. In the order conditions, in each sequence of matrix multiplies, the color indices are compatible according to matrix multiplication rules.
\end{theorem}
The GARK-ROW order four conditions read:
\begin{subequations}
\label{eqn:GARK-ROW-order-conditions}
\begin{alignat}{3}
\label{eqn:GARK-ROW-order1-conditions}
\quad\textnormal{   order 1:   } &
\begin{cases} \mathbf{b}^{\{m\}\,T} \, \one^{\{m\}} = 1, & \textnormal{for  } m = 1,\dots,\nparts; \end{cases} \\
\label{eqn:GARK-ROW-order2-conditions}
\textnormal{   order 2:   } &
\begin{cases}
\mathbf{b}^{\{m\}\,T} \, \mathbf{c}^{\{m,n\}}  = \sfrac{1}{2}, & \\
\mathbf{b}^{\{m\}\,T} \, \mathbf{g}^{\{m,n\}}  = 0,
& \textnormal{for  } m,n = 1,\dots,\nparts;
\end{cases}  \\
\label{eqn:GARK-ROW-order3-conditions}
\textnormal{   order 3:   } &
\begin{cases}
\mathbf{b}^{\{m\}\,T} \, \bigl( \mathbf{c}^{\{m,n\}} \times \mathbf{c}^{\{m,p\}} \bigr) = \sfrac{1}{3}, &
\mathbf{b}^{\{m\}\,T} \, \boldsymbol{\alpha}^{\{m,n\}}\, \mathbf{c}^{\{n,p\}} = \sfrac{1}{6}, \\
\mathbf{b}^{\{m\}\,T} \, \boldsymbol{\gamma}^{\{m,n\}}\,\mathbf{c}^{\{n,p\}} = 0, &
\mathbf{b}^{\{m\}\,T} \, \boldsymbol{\alpha}^{\{m,n\}}\,\mathbf{g}^{\{n,p\}}= 0, \\
\mathbf{b}^{\{m\}\,T} \, \boldsymbol{\gamma}^{\{m,n\}}\, \mathbf{g}^{\{n,p\}} = 0, &
\textnormal{for  } m,n,p = 1,\dots,\nparts;
\end{cases}
\end{alignat}
\begin{gather}
\label{eqn:GARK-ROW-order4-conditions}
\textnormal{order 4:   } \\
\nonumber
\begin{cases}
\mathbf{b}^{\{m\}\,T} \, \left( \mathbf{c}^{\{m,n\}} \times \mathbf{c}^{\{m,p\}} \times \mathbf{c}^{\{m,q\}}\right) = \sfrac{1}{4}, 
&\mathbf{b}^{\{m\}\,T} \, \big((\boldsymbol{\alpha}^{\{m,n\}}\,\mathbf{c}^{\{n,p\}}) \times \mathbf{c}^{\{m,q\}} \big) = \sfrac{1}{8},  \\
\mathbf{b}^{\{m\}\,T} \, \boldsymbol{\alpha}^{\{m,n\}}\, (\mathbf{c}^{\{n,p\}} \times \mathbf{c}^{\{n,q\}}) = \sfrac{1}{12}, 
&\mathbf{b}^{\{m\}\,T}\,\boldsymbol{\alpha}^{\{m,n\}}\,\boldsymbol{\alpha}^{\{n,p\}}\, \mathbf{c}^{\{p,q\}} = \sfrac{1}{24}, \\ 
\mathbf{b}^{\{m\}\,T} \, \big((\boldsymbol{\alpha}^{\{m,n\}}\,\mathbf{g}^{\{n,p\}}) \times \mathbf{c}^{\{m,q\}} \big) = 0, 
&\mathbf{b}^{\{m\}\,T}\, \boldsymbol{\gamma}^{\{m,n\}}\, (\mathbf{c}^{\{n,p\}} \times \mathbf{c}^{\{n,q\}}) = 0, \\ 
\mathbf{b}^{\{m\}\,T} \,\boldsymbol{\gamma}^{\{m,n\}}\,\boldsymbol{\alpha}^{\{n,p\}}\, \mathbf{c}^{\{p,q\}} = 0, 
&\mathbf{b}^{\{m\}\,T}\,\boldsymbol{\alpha}^{\{m,n\}}\,\boldsymbol{\gamma}^{\{n,p\}}\,\mathbf{c}^{\{p,q\}} = 0, \\ 
\mathbf{b}^{\{m\}\,T}\,\boldsymbol{\alpha}^{\{m,n\}}\,\boldsymbol{\alpha}^{\{n,p\}}\,\mathbf{g}^{\{p,q\}} = 0, 
&\mathbf{b}^{\{m\}\,T}\,\boldsymbol{\gamma}^{\{m,n\}}\,\boldsymbol{\alpha}^{\{n,p\}}\,\mathbf{g}^{\{p,q\}} = 0, \\ 
\mathbf{b}^{\{m\}\,T}\,\boldsymbol{\alpha}^{\{m,n\}}\,\boldsymbol{\gamma}^{\{n,p\}}\, \mathbf{g}^{\{p,q\}} = 0,  
&\mathbf{b}^{\{m\}\,T}\,\boldsymbol{\gamma}^{\{m,n\}}\,\boldsymbol{\gamma}^{\{n,p\}}\,\mathbf{c}^{\{p,q\}} = 0, \\ 
\mathbf{b}^{\{m\}\,T}\,\boldsymbol{\gamma}^{\{m,n\}}\,\boldsymbol{\gamma}^{\{n,p\}}\, \mathbf{g}^{\{p,q\}} = 0, &
\textnormal{for  } m,n,p,q = 1,\dots,\nparts.
\end{cases}
\end{gather}
%
%
\end{subequations}
%

\subsection{Internal consistency}

\begin{definition}[Internal consistency]
A partitioned ROW method is internally consistent if:
\begin{subequations}
\label{eqn:ROW-internal-consistency}
\begin{alignat}{4}
\label{eqn:ROW-internal-consistency-alpha}
\mathbf{c}^{\{m,n\}} &= \boldsymbol{\alpha}^{\{m,n\}} \, \one^{\{n\}} = \mathbf{c}^{\{m\}}, \qquad 
& \textnormal{for } m,n& =1,\dots,\nparts, \\
\label{eqn:ROW-internal-consistency-gamma}
\mathbf{g}^{\{m,n\}} &= \boldsymbol{\gamma}^{\{m,n\}} \, \one^{\{n\}} = \mathbf{g}^{\{m\}}, \qquad 
& \textnormal{for } m,n& =1,\dots,\nparts.
\end{alignat}
\end{subequations}
\end{definition}
The order conditions simplify considerably for internally consistent partitioned ROW methods. 

Consider a non-autonomous additively partitioned system \cref{eqn:additively-partitioned-ode} where each component $\fun^{\{m\}}(t,\yy)$ depends explicitly on time. Transform it to autonomous form by adding $t$ to the state, and appending the additively partitioned equation for the time variable $t' = \sum_{m=1}^\nparts \tau^{\{m\}} = 1$.
The stage computation of the GARK-ROS method \cref{eqn:GARK-ROS-stage} applied to non-autonomous system \cref{eqn:additively-partitioned-ode} reads:
\begin{equation}
\begin{split}
\mathbf{k}^{\{q\}} &= h\,\fun^{\{q\}}\mleft( 
\one^{\{q\}} \, t_{n}  + h\,\sum_{m=1}^\nparts \mathbf{c}^{\{q,m\}} \, \tau^{\{m\}},
\one^{\{q\}} \otimes \y_{n}  + \sum_{m=1}^\nparts \boldsymbol{\alpha}^{\{q,m\}} \kron{\nvar}  \, \mathbf{k}^{\{m\}} \mright)    \\
\label{eqn:ROS-non-autonomous}
& \quad +  \bigl( \Id_{s^{\{q\}}} \otimes h\,\Jac_n^{\{q\}} \bigr)\,\sum_{m=1}^\nparts  \boldsymbol{\gamma}^{\{q,m\}} \kron{\nvar}\, \mathbf{k}^{\{m\}} \\
& \quad +  \bigl( \one^{\{q\}} \otimes h^2\,\fun^{\{q\}}_t(t_n,\yy_n) \bigr)\,\sum_{m=1}^\nparts  \mathbf{g}^{\{q,m\}} \tau^{\{m\}}, \qquad q = 1,\dots,\nparts.
\end{split}
\end{equation}
If the internal consistency equation \cref{eqn:ROW-internal-consistency-alpha} holds then the time argument of each function evaluation is $\one^{\{q\}} \, t_{n}  + h\,\mathbf{c}^{\{q\}}$ and is independent of the (arbitrary) partitioning of the time equation. Similarly, if the internal consistency equation \cref{eqn:ROW-internal-consistency-gamma} holds then the coefficient of the time derivative in the stage equation is $\mathbf{g}^{\{q\}}$ and is independent of the partitioning of the time equation.

\begin{remark}[Non-autonomous formulation]
For non-autonomous systems the GARK-ROS method \cref{eqn:GARK-ROS}  computes the set of stages $\mathbf{k}^{\{q\}}$ for process $q$ using the formulation  \cref{eqn:ROS-non-autonomous} with $\tau^{\{q\}} = 1$ and $\tau^{\{m\}} = 0$ for $m \ne q$. This is equivalent with considering a separate time variable for each process. The time argument of each function evaluation in \cref{eqn:ROS-non-autonomous} is $\one^{\{q\}} \, t_{n}  + h\,\mathbf{c}^{\{q,q\}}$, and the coefficient of the time derivative in the stage equation is $\mathbf{g}^{\{q,q\}}$. The same holds for GARK-ROW methods \cref{eqn:GARK-ROW} on non-autonomous systems.
\end{remark}


%
\begin{table}
\centering
\begin{tabular}{|| C{3em} || C{8em}  || C{7em} || C{12em}   ||}
\hline\hline
$\tree$ &  Labels & $F(\tree)$ & $\phi(\tree) \in \{ 1/\gamma(\tree), 0\}$  \\
\hline\hline
$\tree^{\langle w,1 \rangle}_1$ &
\begin{tikzpicture}[grow=up, level distance=0.8cm, scale=0.75]\Tree [.\node[circle,draw,thick,opacity=1,fill=none,inner sep=1.5,minimum size=4,label=right:$\mathbf{b}^{\{m\}T}$](b){m};  ]\end{tikzpicture} &
$\fun^{\{m\}}$ & $\mathbf{b}^{\{m\}T} \, \one^{\{m\}} = 1$  \\ 
\hline\hline
$\tree^{\langle w,1 \rangle}_2 $ &
\begin{tikzpicture}[grow=up, level distance=0.8cm, scale=0.75]\Tree [.\node[circle,draw,thick,opacity=1,fill=none,inner sep=1.5,minimum size=4,label=right:$\mathbf{b}^{\{m\}T}$](j){m}; [.\node[circle,draw,thick,opacity=1,fill=none,inner sep=1.5,minimum size=4,label=right:$\boldsymbol{\alpha}^{\{m,n\}}$](k){n};  ] ]\end{tikzpicture} &    
$\fun_\y^{\{m\}}\,\fun^{\{n\}}$ & $\mathbf{b}^{\{m\}T} \, \boldsymbol{\alpha}^{\{m,n\}} \, \one^{\{n\}} = \frac{1}{2}$  \\ 
$\tree^{\langle w,2 \rangle}_2 $ &
\begin{tikzpicture}[grow=up, level distance=0.8cm, scale=0.75]\Tree [.\node[rectangle,draw,thick,opacity=1,fill=none,inner sep=1.5,minimum size=4,label=right:$\mathbf{b}^{\{m\}T}$](j){m}; [.\node[circle,draw,thick,opacity=1,fill=none,inner sep=1.5,minimum size=4,label=right:$\boldsymbol{\gamma}^{\{m,n\}}$](k){n};  ] ]\end{tikzpicture} &    
$\Lb^{\{m\}}\,\fun^{\{n\}}$ & $\mathbf{b}^{\{m\}T} \, \boldsymbol{\gamma}^{\{m,n\}} \, \one^{\{n\}} = 0$  \\ 
\hline\hline
%
%
%
$\tree^{\langle w,1 \rangle}_{3,1}$ &
\begin{tikzpicture}[grow=up, level distance=0.8cm, scale=0.75]\Tree [.\node[circle,draw,thick,opacity=1,fill=none,inner sep=1.5,minimum size=4,label=right:$\mathbf{b}^{\{m\}T}$](j){m}; [.\node[circle,draw,thick,opacity=1,fill=none,inner sep=1.5,minimum size=4,label=right:$\boldsymbol{\alpha}^{\{m.n\}}$](l){n};  ][.\node[circle,draw,thick,opacity=1,fill=none,inner sep=1.5,minimum size=4,label=right:$\boldsymbol{\alpha}^{\{m,p\}}$](k){p};  ] ]\end{tikzpicture} &   
$\fun_{\y,\y}^{\{m\}}\,(\fun^{\{n\}},\fun^{\{p\}}) $ & $\begin{array}{l} \mathbf{b}^{\{m\}T} \, \big( (\boldsymbol{\alpha}^{\{m,n\}}\, \one^{\{n\}}) \\ \times (\boldsymbol{\alpha}^{\{m,p\}}\, \one^{\{p\}}) \big) = \frac{1}{3} \end{array}$  \\ 
\hline
$\tree^{\langle w,1 \rangle}_{3,2}$ &
\begin{tikzpicture}[grow=up, level distance=0.8cm, scale=0.75]\Tree [.\node[circle,draw,thick,opacity=1,fill=none,inner sep=1.5,minimum size=4,label=right:$\mathbf{b}^{\{m\}T}$](j){m}; [.\node[circle,draw,thick,opacity=1,fill=none,inner sep=1.5,minimum size=4,label=right:$\boldsymbol{\alpha}^{\{m,n\}}$](k){n}; [.\node[circle,draw,thick,opacity=1,fill=none,inner sep=1.5,minimum size=4,label=right:$\boldsymbol{\alpha}^{\{n,p\}}$](l){p};  ] ] ]\end{tikzpicture} &
$\fun_\y^{\{m\}}\,\fun_\y^{\{n\}}\,\fun^{\{p\}}$ & $\begin{array}{l} \mathbf{b}^{\{m\}T} \, \boldsymbol{\alpha}^{\{m,n\}}\\ \cdot \boldsymbol{\alpha}^{\{n,p\}}\,\one^{\{p\}} = \frac{1}{6} \end{array}$  \\ 
\hline
$\tree^{\langle w,2 \rangle}_{3,2}$ &
\begin{tikzpicture}[grow=up, level distance=0.8cm, scale=0.75]\Tree [.\node[circle,draw,thick,opacity=1,fill=none,inner sep=1.5,minimum size=4,label=right:$\mathbf{b}^{\{m\}T}$](j){m}; [.\node[rectangle,draw,thick,opacity=1,fill=none,inner sep=1.5,minimum size=4,label=right:$\boldsymbol{\alpha}^{\{m,n\}}$](k){n}; [.\node[circle,draw,thick,opacity=1,fill=none,inner sep=1.5,minimum size=4,label=right:$\boldsymbol{\gamma}^{\{n,p\}}$](l){p};  ] ] ]\end{tikzpicture} &
$\fun_\y^{\{m\}}\,\Lb^{\{n\}}\,\fun^{\{p\}}$ & $\begin{array}{l} \mathbf{b}^{\{m\}T} \, \boldsymbol{\alpha}^{\{m,n\}}\\ \cdot \boldsymbol{\gamma}^{\{n,p\}}\,\one^{\{p\}} = 0 \end{array}$  \\ 
\hline
$\tree^{\langle w,3 \rangle}_{3,2}$ &
\begin{tikzpicture}[grow=up, level distance=0.8cm, scale=0.75]\Tree [.\node[rectangle,draw,thick,opacity=1,fill=none,inner sep=1.5,minimum size=4,label=right:$\mathbf{b}^{\{m\}T}$](j){m}; [.\node[circle,draw,thick,opacity=1,fill=none,inner sep=1.5,minimum size=4,label=right:$\boldsymbol{\gamma}^{\{m,n\}}$](k){n}; [.\node[circle,draw,thick,opacity=1,fill=none,inner sep=1.5,minimum size=4,label=right:$\boldsymbol{\alpha}^{\{n,p\}}$](l){p};  ] ] ]\end{tikzpicture} &
$\Lb^{\{m\}}\,\fun_\y^{\{n\}}\,\fun^{\{p\}}$ & $\begin{array}{l} \mathbf{b}^{\{m\}T} \, \boldsymbol{\gamma}^{\{m,n\}}\\ \cdot \boldsymbol{\alpha}^{\{n,p\}}\,\one^{\{p\}} = 0 \end{array}$  \\ 
\hline
$\tree^{\langle w,4 \rangle}_{3,2}$ &
\begin{tikzpicture}[grow=up, level distance=0.8cm, scale=0.75]\Tree [.\node[rectangle,draw,thick,opacity=1,fill=none,inner sep=1.5,minimum size=4,label=right:$\mathbf{b}^{\{m\}T}$](j){m}; [.\node[rectangle,draw,thick,opacity=1,fill=none,inner sep=1.5,minimum size=4,label=right:$\boldsymbol{\gamma}^{\{m,n\}}$](k){n}; [.\node[circle,draw,thick,opacity=1,fill=none,inner sep=1.5,minimum size=4,label=right:$\boldsymbol{\gamma}^{\{n,p\}}$](l){p};  ] ] ]\end{tikzpicture} &
$\Lb^{\{m\}}\,\Lb^{\{n\}}\,\fun^{\{p\}}$ & $\begin{array}{l} \mathbf{b}^{\{m\}T} \, \boldsymbol{\gamma}^{\{m,n\}}\\ \cdot \boldsymbol{\gamma}^{\{n,p\}}\,\one^{\{p\}} = 0 \end{array}$  \\ 
\hline\hline
\end{tabular}
\caption{$\NTW$ trees of orders 1 to 3 for the GARK-ROW numerical solution. Square vertices of color $\nu$ correspond to $\Lb^{\{\nu\}}$, and round vertices to derivatives of $\fun^{\{\nu\}}$. A root of color $m$ is labelled $\mathbf{b}^{\{m\}T}$. Nodes of color $m$ with a round parent of color $q$ are labelled $\boldsymbol{\alpha}^{\{q,m\}}$. Nodes of color $m$ with a square parent of color $q$ are labelled $\boldsymbol{\gamma}^{\{q,m\}}$.}
\label{tab:GARK-ROW-trees_1-3_numeric}
\end{table}


\section{Linear stability}
\label{sec:linear-stability}

Consider the scalar test problem
\begin{equation}
\label{eqn:scalar-dahlquist}
\y' = \lambda^{\{1\}}\,\y + \cdots + \lambda^{\{\nparts\}}\, \y.
\end{equation}
Application of the GARK-ROS method \cref{eqn:GARK-ROW} to \cref{eqn:scalar-dahlquist} leads to the same stability equation as the application of a GARK scheme.
%
Using the notation \cref{eqn:GARK-Rosenbrock-butcher} and defining $\mathbf{B} = \mathbf{A} + \mathbf{G} \in \Re^{s \times s}$, and
\begin{equation*}
\begin{split}
z^{\{m\}} &\coloneqq h\,\lambda^{\{m\}}, \quad s \coloneqq \sum_{m=1}^\nparts s^{\{m\}}, \quad
Z \coloneqq \operatorname{diag}_{m=1,\dots,\nparts} \, \big\{ z^{\{m\}}\,\Id_{s^{\{m\}}} \big\} \in \Re^{s \times s},
\end{split}
\end{equation*}
we obtain $\y_{n+1} = R(Z)\,\y_{n}$, with
\begin{equation}
\label{eqn:GARK-ROS-stability}
R(Z) = 1 + \b^{T} \, \left( \Id_{s} - Z\,\mathbf{B} \right)^{-1}\,Z\,\one_{s}
 = 1 + \b^{T} \,Z\, \left( \Id_{s} - \mathbf{B}\,Z \right)^{-1}\,\one_{s},
\end{equation}
which equals the stability function of a GARK scheme with coefficients $(\mathbf{B},\b)$. The following definition extends immediately from GARK to GARK-ROS schemes.
\begin{definition}[Stiff accuracy]
Let $\mathbf{e}_{s} \in \Re^s$ be a vector with the last entry equal to one, and all other entries equal to zero.
The GARK-ROS method \cref{eqn:GARK-ROW} is called stiffly accurate if
\begin{equation*}
\b^{T} = \mathbf{e}_{s}^{T}\,\mathbf{B}
\quad \Leftrightarrow \quad
\b^{\{q\}T} = \mathbf{e}_{s^{\{\nparts\}}}^{T}\,\boldsymbol{\beta}^{\{\nparts,q\}}, \quad q = 1,\dots,\nparts.
\end{equation*}
\end{definition}
For a stiffly accurate GARK-ROS scheme the stability function \cref{eqn:GARK-ROS-stability} becomes:
\begin{equation}
\label{eqn:GARK-ROS-stability-SA}
\begin{split}
R(Z)  = z_{\nparts}^{-1}\,\mathbf{e}_{s}^{T}\,\left( Z^{-1} - \mathbf{B} \right)^{-1}\,\one_{s}.
\end{split}
\end{equation}
If 
$\diag(1/z_1,\ldots,1/z_{n-1},0)-\mathbf{B}$ is nonsingular then $R(Z) \to 0$ when $z_{\nparts} \to \infty$. This condition is automatically fulfilled for decoupled GARK-ROW schemes discussed in \cref{sec:decoupled-GARK-ROW}.

\section{GARK-ROW multimethods}
\label{sec:multimethods-GARK-ROW}

The GARK-ROS/GARK-ROW framework allows to construct different types of multimethods. In the following we address decoupled and  linearly implicit-explicit (for short, LIMEX) GARK-ROW schemes, as well as implicit/linearly implicit GARK methods arising from the GARK-ROS framework.

\subsection{Decoupled IMIM GARK-ROW schemes}
\label{sec:decoupled-GARK-ROW}
Consider now an $\nparts$-way additively partitioned system \cref{eqn:additively-partitioned-ode}. Application of a traditional ROW scheme solves a single system with matrix $\Idvar - h\gamma (\Lb^{\{1\}}+\dots+\Lb^{\{\nparts\}})$. The GARK-ROW scheme \cref{eqn:GARK-ROW} applied to the $\nparts$-way partitioned system reads \cref{eqn:GARK-ROW-stage}: 
\begin{equation}
\label{eqn:GARK-ROW-FULL}
\begin{split}
& \left( \Id_{s\nvar}  - \diag_q \bigl\{ \Id_{s^{\{q\}}} \otimes h\,\Lb^{\{q\}} \bigr\} \cdot \left( \mathbf{G} \otimes \Idvar \right) \right) \cdot \mathbf{K} 
= h\,\mathbf{F}\Bigl(  \one_s \otimes \y_{n}  +   \mathbf{A} \kron{\nvar} \mathbf{K} \Bigr),  \\
& \y_{n+1} = \y_{n} + \b^{T} \kron{\nvar}  \mathbf{K}, \\
& \mathbf{K} \coloneqq \begin{bmatrix} \mathbf{k}^{\{1\}} \\ \smash[t]{\vdots} \\ \mathbf{k}^{\{\nparts\}} \end{bmatrix}, \quad  
\mathbf{F} \coloneqq \begin{bmatrix} \fun^{\{1\}}\mleft( \one^{\{1\}} \otimes \y_{n}  + \mathbf{A}^{\{1,:\}} \kron{\nvar} \mathbf{K} \mright) \\ 
\smash[t]{\vdots} \\ 
\fun^{\{\nparts\}}\mleft( \one^{\{\nparts\}} \otimes \y_{n}  + \mathbf{A}^{\{\nparts,:\}} \kron{\nvar} \mathbf{K} \mright) \end{bmatrix}.
\end{split}
\end{equation}
Equation \cref{eqn:GARK-ROW-FULL} shows that the stage vectors are obtained by solving a linear system of dimension $sd$ that, in general, couples all components together.
%
To increase computational efficiency we look for schemes where each stage $\mathbf{k}_i^{\{q\}}$ is obtained by solving a $\nvar$-dimensional linear system with matrix $\Id_{\nvar} - h\,  \gamma_{i,i}^{\{q,q\}}\,\Lb^{\{q\}}$. We call such methods ``decoupled.''


\begin{definition}[Decoupled schemes]
GARK-ROW schemes \cref{eqn:GARK-ROW} are decoupled if they solve the stage equations implicitly in either one process or the other, but not in both in the same time.
\end{definition}

\begin{theorem}
\label{thm:stage-reordering-decoupled}
A method \cref{eqn:GARK-ROW} is decoupled iff there is a permutation vector $v$ (representing the order of stage evaluations), and an associated permutation matrix $\mathcal{V}$,  such that the matrices of coefficients \cref{eqn:GARK-Rosenbrock-butcher} with reordered rows and columns have the following structure:  $\mathcal{V}$ is strictly lower triangular and $\mathbf{G}(v,v)=\mathcal{V}\,\mathbf{G}\,\mathcal{V}$ is lower triangular. 
\end{theorem}
\begin{proof}
Stage reordering $\mathbf{K} \to \mathcal{V}\kron{\nvar} \mathbf{K}$ leads to linear systems \cref{eqn:GARK-ROW-FULL} that are block lower triangular; the argument of the right hand side function also involves a block strictly lower triangular matrix of coefficients, and therefore \cref{eqn:GARK-ROW-FULL} can be solved by forward substitution.
\end{proof}

To illustrate how the property in \cref{thm:stage-reordering-decoupled} applies, consider the scalar formulation of the stage computations \cref{eqn:GARK-ROW-stage}:
\begin{align*}
k_i^{\{q\}} &= h\,\fun^{\{q\}} \mleft(\y_{n} +  \sum_{m=1}^{q-1}  \sum_{j=1}^{i} \alpha_{i,j}^{\{q,m\}}\,k_j^{\{m\}} +   \sum_{m=q}^\nparts  \sum_{j=1}^{i-1 } \alpha_{i,j}^{\{q,m\}}\,k_j^{\{m\}} \mright) \\
\nonumber
& \quad + h \, \Lb^{\{q\}}\,\sum_{m = 1}^q \sum_{j=1}^i \gamma_{i,j}^{\{q,m\}}\,k_j^{\{m\}}
      + h \, \Lb^{\{q\}}\,\sum_{m = q+1}^\nparts \sum_{j=1}^{i-1} \gamma_{i,j}^{\{q,m\}}\,k_j^{\{m\}}.
\end{align*}
The stages are solved in the order $k_i^{\{1\}},\dots,k_i^{\{\nparts\}}$, then $k_{i+1}^{\{1\}},\dots,k_{i+1}^{\{\nparts\}}$, etc. The computation of stage $k_i^{\{q\}}$ uses all $k_j^{\{1\}},\dots,k_j^{\{\nparts\}}$ for $j < i$, as well as $k_i^{\{1\}},\dots,k_i^{\{q-1\}}$, which have already been computed. Stage $k_i^{\{q\}}$ is obtained by solving a linear system with matrix $\Id_{s^{\{q\}}} - h\,  \gamma_{i,i}^{\{q,q\}}\,\Lb^{\{q\}}$.
Here we allow $\boldsymbol{\alpha}^{\{q,m\}}$ for $m<q$ to be lower triangular and do not demand a strictly lower triangular structure.

\begin{remark}
If the coefficient matrices $\boldsymbol{\gamma}^{\{\ell,m\}}$  are strictly lower triangular for all $\ell \ne m$ then all implicit stages $\mathbf{k}^{\{\ell\}}_i$, $\ell=1,\dots,\nparts$, can be evaluated in parallel.
\end{remark}

\begin{remark}[First special case]
A first interesting special case arises when:
\begin{equation*}
\begin{split}
\boldsymbol{\gamma}^{\{q,m\}} &= \begin{cases}
\underline{\boldsymbol{\gamma}}~ \textnormal{(lower triangular)},  & m = 1, \dots, q-1, \\
\boldsymbol{\gamma}~ \textnormal{(lower triangular)},  & m = q, \\
\overline{\boldsymbol{\gamma}}~ \textnormal{(strictly lower triangular)}, & m = q+1, \dots, \nparts,
\end{cases}
\\
\boldsymbol{\alpha}^{\{q,m\}}  &=  \begin{cases}
\underline{\boldsymbol{\alpha}}~ \textnormal{(lower triangular)},  & m = 1, \dots, q-1, \\
\boldsymbol{\alpha}~ \textnormal{(strictly lower triangular)},  & m = q, \\
\overline{\boldsymbol{\alpha}}~ \textnormal{(strictly lower triangular)}, & m = q+1, \dots, \nparts.
\end{cases}
\end{split}
\end{equation*}
The computations are carried out as follows:
\begin{align*}
k_i^{\{q\}} &= h\,\fun^{\{q\}} \mleft(\y_{n} +   \sum_{j=1}^{i} \underline{\alpha}_{i,j}\,\biggl(\sum_{m=1}^{q-1} k_j^{\{m\}}\biggr) + \sum_{j=1}^{i-1} \alpha_{i,j} \,k_j^{\{q\}}+ \sum_{j=1}^{i-1 } \overline{\alpha}_{i,j} \, \biggl( \sum_{m=q+1}^\nparts  k_j^{\{m\}}\biggr) \mright) \\
\nonumber
& \quad + h \, \Lb^{\{q\}}\,\left( \sum_{j=1}^i \underline{\gamma}_{i,j}\,\biggl(\sum_{m=1}^{q-1} k_j^{\{m\}}\biggr)+ \sum_{j=1}^i \gamma_{i,j}\,k_j^{\{q\}}
      + \sum_{j=1}^{i-1} \overline{\gamma}_{i,j}\, \biggl( \sum_{m=q+1}^\nparts  k_j^{\{m\}}\biggr) \right).
\end{align*}
\ifreport
The Butcher tableau \cref{eqn:GARK-Rosenbrock-butcher} reads:
\begin{equation}
\begin{butchertableau}{c|c}
 \A  & \G \\ 
\hline
\b^T &
\end{butchertableau}
~=
~\raisebox{25.5pt}{$
\begin{butchertableau}{cccc|cccc}
\boldsymbol{\alpha} & \overline{\boldsymbol{\alpha}} & \ldots & \overline{\boldsymbol{\alpha}} & \boldsymbol{\gamma} & \overline{\boldsymbol{\gamma}} & \ldots & \overline{\boldsymbol{\gamma}} \\
\underline{\boldsymbol{\alpha}} &\boldsymbol{\alpha} & \ldots & \overline{\boldsymbol{\alpha}}  & \underline{\boldsymbol{\gamma}} &\boldsymbol{\gamma} & \ldots & \overline{\boldsymbol{\gamma}} \\
\vdots & \vdots & & \vdots  & \vdots & \vdots & & \vdots  \\
\underline{\boldsymbol{\alpha}} & \underline{\boldsymbol{\alpha}} & \ldots & \boldsymbol{\alpha}  & \underline{\boldsymbol{\gamma}} & \underline{\boldsymbol{\gamma}} & \ldots & \boldsymbol{\gamma} \\ 
\hline
\b^{\{1\}T} & \b^{\{2\}T} & \ldots &\b^{\{N\}T}
\end{butchertableau}
$}.
\end{equation}
\fi
\end{remark}

\begin{remark}[Second special case]
\label{rem:special-case-2}
A second interesting case arises when:
\begin{subequations}
\label{eqn:special-case-2}
\begin{equation}
\underline{\boldsymbol{\alpha}} = \boldsymbol{\overline{\alpha}}, \quad 
\underline{\boldsymbol{\gamma}} = \boldsymbol{\overline{\gamma}}~~ \textnormal{(strictly lower triangular)}; \quad
\boldsymbol{b}^{\{q\}} =  \boldsymbol{b} ~~ \forall\,q, \quad
\boldsymbol{\overline{c}} =  \boldsymbol{c}.
\end{equation}
The computations are carried out as follows:
\begin{align}
k_i^{\{q\}} &= h\,\fun^{\{q\}} \mleft(\y_{n} +  \sum_{j=1}^{i-1} \alpha_{i,j}\,k_j^{\{q\}} +   \sum_{j=1}^{i-1 } \overline{\alpha}_{i,j} \,
\biggl( \sum_{m \ne q}  k_j^{\{m\}} \biggr) \mright) \\
\nonumber
& \quad + h \, \Lb^{\{q\}}\,\left( \sum_{j=1}^i \gamma_{i,j}\,k_j^{\{q\}}
      +\sum_{j=1}^{i-1} \overline{\gamma}_{i,j}\,\biggl( \sum_{m \ne q}  k_j^{\{m\}} \biggr) \right), \quad q = 1,\dots,\nparts, \\
\y_{n+1} &= \y_{n} +  \sum_{i=1}^{s} b_i^T \, \biggl( \sum_{m=1}^{\nparts}  k_i^{\{m\}} \biggr).   
\end{align}
\end{subequations}
Here $(\boldsymbol{b},\boldsymbol{\alpha},\boldsymbol{\gamma})$ is a base Rosenbrock or Rosenbrock-W scheme, and $(\boldsymbol{b},\boldsymbol{\overline{\alpha}},\boldsymbol{\overline{\gamma}})$ are the coupling coefficients.
\end{remark}
The GARK-ROW order three conditions \cref{eqn:GARK-ROW-order3-conditions} for methods \cref{eqn:special-case-2} are as follows. Both the base scheme $(\boldsymbol{b},\boldsymbol{\alpha},\boldsymbol{\gamma})$ and coupling scheme $(\boldsymbol{b},\boldsymbol{\overline{\alpha}},\boldsymbol{\overline{\gamma}})$ need to be order three  Rosenbrock-W schemes. The following third order coupling conditions are also needed:
\begin{equation}
\label{eqn:ROW-coupling-order-3}
\mathbf{b}^T \, \boldsymbol{\alpha}\,\overline{\mathbf{g}}= 
\mathbf{b}^T \, \overline{\boldsymbol{\alpha}}\,\mathbf{g}= 
\mathbf{b}^T \, \boldsymbol{\gamma}\, \overline{\mathbf{g}} = 
\mathbf{b}^T \, \overline{\boldsymbol{\gamma}}\, \mathbf{g} = 0.
\end{equation}
Choosing $\overline{\boldsymbol{\gamma}} = \boldsymbol{0}$ means that $(\boldsymbol{b},\boldsymbol{\overline{\alpha}})$ is an explicit Runge--Kutta scheme, and only the coupling equation $\mathbf{b}^T \, \overline{\boldsymbol{\alpha}}\,\mathbf{g}=0$ needs to be imposed. 

The GARK-ROS order four conditions \cref{eqn:GARK-ROS-order-conditions} for methods \cref{eqn:special-case-2} require that the base and the coupling schemes are order four Rosenbrock methods.  In addition, one needs to satisfy the third order coupling conditions:
\begin{equation}
\label{eqn:ROS-coupling-order-3}
\mathbf{b}^{T} \, \boldsymbol{\beta}\, \overline{\mathbf{e}} =
\mathbf{b}^{T} \, \overline{\boldsymbol{\beta}}\, \mathbf{e} = \sfrac{1}{6},
\end{equation}
as well as the fourth order coupling conditions: 
\begin{equation}
\label{eqn:ROS-coupling-order-4}
\begin{split}
&\mathbf{b}^{T} \, \big(({\boldsymbol{\alpha}}\,\overline{\mathbf{e}}) \times  \mathbf{c} \big) = 
\mathbf{b}^{T} \, \big((\overline{\boldsymbol{\alpha}}\,\mathbf{e}) \times  \mathbf{c} \big)  =\sfrac{1}{8}, \\
&\mathbf{b}^{T} \,\overline{\boldsymbol{\beta}}\,\boldsymbol{\beta}\,\mathbf{e} = 
\mathbf{b}^{T} \,\overline{\boldsymbol{\beta}}\,\overline{\boldsymbol{\beta}}\,\mathbf{e} =\mathbf{b}^{T} \, \overline{\boldsymbol{\beta}}\,\boldsymbol{\beta}\,\overline{\mathbf{e}} = 
\mathbf{b}^{T} \,\boldsymbol{\beta}\,\overline{\boldsymbol{\beta}}\,\overline{\mathbf{e}} = \mathbf{b}^{T} \,\boldsymbol{\beta}\,\overline{\boldsymbol{\beta}}\,\mathbf{e} =\mathbf{b}^{T} \,\boldsymbol{\beta}\,\boldsymbol{\beta}\,\overline{\mathbf{e}} = \sfrac{1}{24}.
\end{split}
\end{equation} 
Choosing $\overline{\boldsymbol{\gamma}} = \boldsymbol{0}$ further simplifies the coupling equations \cref{eqn:ROS-coupling-order-3} and \cref{eqn:ROS-coupling-order-4}.

For decoupled GARK-ROS/ROW schemes the stability function \cref{eqn:GARK-ROS-stability} is rewritten using the permutation matrix from \cref{thm:stage-reordering-decoupled}:
\begin{equation}
\label{eqn:GARK-ROS-decoupled-stability}
\begin{split}
R(Z)  
& = 1 + (\mathcal{V}\,\b)^{T} \, (\mathcal{V}\,Z\,\mathcal{V})\, \left( \mathbf{I}_{s} - (\mathcal{V}\,\mathbf{B}\,\mathcal{V})\,(\mathcal{V}\,Z\,\mathcal{V}) \right)^{-1}\,\one_{s}.
\end{split}
\end{equation}
The matrix $\mathcal{V}\,\mathbf{B}\,\mathcal{V}$ is lower triangular, with the diagonal entries equal to the diagonal entries of $\G^{\{m,m\}}$. 
The stability function \cref{eqn:GARK-ROS-decoupled-stability} is a rational function of the form:
\begin{equation}
\begin{split}
R(Z)   & = \frac{\varphi(z^{\{1\}},\dots,z^{\{\nparts\}})}{\prod_{m=1}^\nparts\, \prod_{i=1}^{s^{\{m\}}}\, (1 - \gamma^{\{m,m\}}_{i,i} \,z^{\{m\}})}.
\end{split}
\end{equation}

\subsection{Implicit-explicit (IMEX) GARK-ROW schemes}
\label{sec:IMEX-GARK-ROW}

Consider now a two-way partitioned system driven by a non-stiff component $\fun\E$ and a stiff component $\fun\I$:
\begin{equation}
\label{eqn:imex-ode}
\y' = \fun\E \left(\y\right) + \fun\I \left(\y\right).
\end{equation}
We consider a GARK-ROW scheme \cref{eqn:GARK-ROW} applied to \cref{eqn:imex-ode} that has the form:
\ifreport
\begin{subequations}
\label{eqn:GARK-ROW-IMEX-scalar}
\begin{align}
\label{eqn:GARK-ROW-IMEX-scalar-stages-E}
k_i\E &= h\,\fun\E \mleft(\y_{n} +  \sum_{j=1}^{i-1} \alpha_{i,j}\EE\,k_j\E
+  \sum_{j=1}^{i-1 } \alpha_{i,j}\EI\,k_j\I \mright) \\
\nonumber & + h \, \Lb\E\,\left( \sum_{j=1}^{i-1} \gamma_{i,j}\EE\,k_j\E + \sum_{j=1}^{i-1} \gamma_{i,j}\EI\,k_j\I \right), \\
\label{eqn:GARK-ROW-IMEX-scalar-stages-I}
k_i\I &= h\,\fun\I \mleft(\y_{n} +   \sum_{j=1}^{i} \alpha_{i,j}\IE\,k_j\E
+  \sum_{j=1}^{i-1} \alpha_{i,j}\II\,k_j\I \mright) \\
\nonumber & \quad + h \, \Lb\I\,\left( \sum_{j=1}^{i} \gamma_{i,j}\IE\,k_j\E + \sum_{j=1}^{i} \gamma_{i,j}\II\,k_j\I \right), \\
\label{eqn:GARK-ROW-IMEX-scalar-solution}
\qquad y_{n+1} & = y_{n} + \sum_{i=1}^{s\E} b_{i}\E\,k_i\E + \sum_{i=1}^{s\I} b_{i}\I\,k_i\I.
\end{align}
\end{subequations}
In matrix notation the IMEX GARK-ROW scheme \cref{eqn:GARK-ROW-IMEX-scalar} reads:
\fi
\begin{subequations}
\label{eqn:GARK-ROW-IMEX}
\begin{align}
\label{eqn:GARK-ROW-IMEX-stage}
\mathbf{k}\E &= h\,\fun\E\mleft( 
\one_s \otimes \y_{n}  + \boldsymbol{\alpha}\EE \kron{\nvar} \mathbf{k}\E 
+  \boldsymbol{\alpha}\EI \kron{\nvar}  \mathbf{k}\I \mright)  \\
& \quad + \bigl( \Idstage \otimes h\,\mathbf{L}\E \bigr)\,
\left( \boldsymbol{\gamma}\EE \kron{\nvar}   \mathbf{k}\E
+  \boldsymbol{\gamma}\EI \kron{\nvar}  \mathbf{k}\I \right), \\ 
\mathbf{k}\I &= h\,\fun\I\mleft( 
\one_s \otimes \y_{n}  + \boldsymbol{\alpha}\IE \kron{\nvar} \mathbf{k}\E 
+ \boldsymbol{\alpha}\II \kron{\nvar} \mathbf{k}\I \mright)  \\
\nonumber
& \quad + \bigl( \Idstage \otimes h\,\mathbf{L}\I \bigr)\,
\left( \boldsymbol{\gamma}\IE \kron{\nvar}   \mathbf{k}\E
+  \boldsymbol{\gamma}\II \kron{\nvar}  \mathbf{k}\I \right), \\ 
\label{eqn:GARK-ROW-IMEX-solution}
\y_{n+1} &= \y_{n} + \b\E\,\!^{T} \kron{\nvar} \mathbf{k}\E
+  \b\I\,\!^{T} \kron{\nvar} \mathbf{k}\I,
\end{align}
\end{subequations} 
with $\boldsymbol{\alpha}\EE$, $\boldsymbol{\alpha}\EI$, $\boldsymbol{\alpha}\II$, $\boldsymbol{\gamma}\EE$, $\boldsymbol{\gamma}\EI$ strictly lower triangular, and $\boldsymbol{\alpha}\IE$, $\boldsymbol{\gamma}\IE$, $\boldsymbol{\gamma}\II$ lower triangular. 

Of particular interest are schemes where $\boldsymbol{\gamma}\EE = \boldsymbol{0}$, $\boldsymbol{\gamma}\EI = \boldsymbol{0}$, which leads to the structure of ERK-ROS methods proposed in \citep{Hai_2014_ERK-ROS}. In this case the non-stiff component $\fun\E$ is solved with an explicit GARK scheme, and the stiff component $\fun\I$ with a linearly implicit scheme.

For order three $(\mathbf{b}\E,\boldsymbol{\alpha}\EE)$ needs to be a third order explicit Runge--Kutta scheme. For arbitrary Jacobian approximations $\Lb\I$ the scheme \linebreak $(\mathbf{b}\I,\boldsymbol{\alpha}\II,\boldsymbol{\gamma}\II)$ has to be a third order Rosenbrock-W method. In addition, assuming the internal consistency  \cref{eqn:ROW-internal-consistency-alpha}, the coupling order three conditions \cref{eqn:GARK-ROW-order3-conditions} are:
\begin{equation}
\label{eqn:IMEX-ROW-third-order-coupling}
\begin{alignedat}{3}
&\mathbf{b}\E\,\!^T \, \boldsymbol{\alpha}\EI\, \mathbf{c}\I = \sfrac{1}{6}, & \qquad &
\mathbf{b}\E\,\!^T \, \boldsymbol{\alpha}\EI\,\mathbf{g}\I = 0, \\
&\mathbf{b}\I\,\!^T \, \boldsymbol{\alpha}\IE\, \mathbf{c}\E = \sfrac{1}{6}, & \qquad &
\mathbf{b}\I\,\!^T \, \boldsymbol{\gamma}\IE\,\mathbf{c}\E = 0.
\end{alignedat}
\end{equation}
If the exact Jacobian is used, $\Lb\I = \Jac_{n}\I$, then the implicit scheme needs to be a third order Rosenbrock method, and the coupling conditions are:
\begin{equation}
\label{eqn:IMEX-ROS-third-order-coupling}
\mathbf{b}\E\,\!^T \, \boldsymbol{\alpha}\EI\, \mathbf{e}\I = \sfrac{1}{6}, \qquad
\mathbf{b}\I\,\!^T \, \boldsymbol{\beta}\IE\, \mathbf{c}\E = \sfrac{1}{6}.
\end{equation}

When the exact Jacobian is used, for order four one needs $(\mathbf{b}\E,\boldsymbol{\alpha}\EE)$ to be a fourth order explicit Runge--Kutta scheme, and $(\mathbf{b}\I,\boldsymbol{\alpha}\II,\boldsymbol{\gamma}\II)$ to be a fourth order Rosenbrock method. In this case
the coupling order four conditions are:
\begin{equation}
\begin{alignedat}{3}
&\mathbf{b}\E\,\!^T \, \big((\boldsymbol{\alpha}\EI\,\mathbf{e}\I) \times \mathbf{c}\E \big) = \sfrac{1}{8}, & \qquad &
\mathbf{b}\I\,\!^T \, \big((\boldsymbol{\alpha}\IE\,\mathbf{c}\E) \times \mathbf{c}\I \big) = \sfrac{1}{8}, \\
&\mathbf{b}\E\,\!^T \, \boldsymbol{\alpha}\EI\, (\mathbf{c}\I)^{\times 2}  = \sfrac{1}{12}, & \qquad & 
\mathbf{b}\I\,\!^T \, \boldsymbol{\beta}\IE\, (\mathbf{c}\E)^{\times 2}  = \sfrac{1}{12}, \\ 
\label{eqn:IMEX-ROS-fourth-order-coupling}
&\mathbf{b}\E\,\!^T\,\boldsymbol{\alpha}\EE\,\boldsymbol{\alpha}\EI\, \mathbf{e}\I = \sfrac{1}{24}, &
&\mathbf{b}\E\,\!^T\,\boldsymbol{\alpha}\EI\,\boldsymbol{\beta}\IE\, \mathbf{c}\E = \sfrac{1}{24}, \\
&\mathbf{b}\E\,\!^T\,\boldsymbol{\alpha}\EI\,\boldsymbol{\beta}\II\, \mathbf{e}\I = \sfrac{1}{24}, &
&\mathbf{b}\I\,\!^T\,\boldsymbol{\beta}\IE\,\boldsymbol{\alpha}\EE\, \mathbf{c}\E = \sfrac{1}{24}, \\
&\mathbf{b}\I\,\!^T\,\boldsymbol{\beta}\IE\,\boldsymbol{\alpha}\EI\, \mathbf{e}\I = \sfrac{1}{24}, &
&\mathbf{b}\I\,\!^T\,\boldsymbol{\beta}\II\,\boldsymbol{\beta}\IE\, \mathbf{c}\E = \sfrac{1}{24}.
\end{alignedat}
\end{equation}

\begin{remark}
\label{rem:imex-special-case}
An interesting special case is when \cref{eqn:GARK-ROW-IMEX} uses:
\begin{equation}
\label{eqn:GARK-ROW-IMEX-case1}
\begin{split}
\boldsymbol{\alpha}\EI &=  \boldsymbol{\alpha}\EE =  \boldsymbol{\alpha}\E, \quad
\boldsymbol{\alpha}\IE = \boldsymbol{\alpha}\II = \boldsymbol{\alpha}\I,\\
\boldsymbol{\gamma}\IE &= \boldsymbol{\gamma}\II = \boldsymbol{\gamma}\I, \quad
\mathbf{g}\I = \boldsymbol{\gamma}\I \, \one, \quad
\mathbf{c} = \boldsymbol{\alpha}\I \, \one = \boldsymbol{\alpha}\E \, \one.
\end{split}
\end{equation}
In this case the method \cref{eqn:GARK-ROW-IMEX} couples an explicit Runge--Kutta scheme
\linebreak $(\mathbf{b}\E,\boldsymbol{\alpha}\E)$ with a Rosenbrock (or Rosenbrock-W) scheme $(\mathbf{b}\I,\boldsymbol{\alpha}\I,\boldsymbol{\gamma}\I)$. 
\ifreport
The computations proceed as follows:
\begin{subequations}
\label{eqn:GARK-ROW-IMEX-computations-case1}
\begin{align}
\label{eqn:GARK-ROW-IMEX-stage-case1}
\mathbf{k}\E &= h\,\fun\E\mleft( 
\one_{s\E} \otimes \y_{n}  + \boldsymbol{\alpha}\E \kron{\nvar}  (\mathbf{k}\E + \mathbf{k}\I)   \mright),  \\
\mathbf{k}\I &= h\,\fun\I\left( 
\one_{s\I} \otimes \y_{n}  +  \boldsymbol{\alpha}\I \kron{\nvar}   (\mathbf{k}\E + \mathbf{k}\I)  \right)  \\
\nonumber
& \quad + \bigl( \Idstage \otimes h\,\mathbf{L}\I \bigr)\,
\left( \boldsymbol{\gamma}\I \kron{\nvar}   (\mathbf{k}\E + \mathbf{k}\I)  \right), \\ 
\label{eqn:GARK-ROW-IMEX-solution-case1}
\y_{n+1} &= \y_{n} + \b\E\,\!^{T} \kron{\nvar}  \mathbf{k}\E
+   \b\I\,\!^{T} \kron{\nvar}  \mathbf{k}\I.
\end{align}
\end{subequations}
\fi
For IMEX order $p$ the explicit and the linearly implicit method need to have order at least $p$.
For arbitrary $\Lb\I$ the $p=3$ GARK-ROW coupling conditions \cref{eqn:GARK-ROW-order3-conditions} are:
\begin{equation}
\label{eqn:GARK-ROW-IMEX-case1-order3}
\mathbf{b}\E\,\!^T \, \boldsymbol{\alpha}\E\,\mathbf{g}\I = 0,
\end{equation}
and the $p=4$ the GARK-ROW coupling conditions  \cref{eqn:GARK-ROW-order4-conditions} simplify to
%
\begin{equation}
\begin{alignedat}{3}
& \mathbf{b}\E\,\!^T \, \big((\boldsymbol{\alpha}\E\,\mathbf{g}\I) \times \mathbf{c} \big) = 0, 
& \qquad & \mathbf{b}\E\,\!^T\,\boldsymbol{\alpha}\E\,\boldsymbol{\alpha}\E\, \mathbf{g}\I = 0, \\
& \mathbf{b}\E\,\!^T\,\boldsymbol{\alpha}\E\,\boldsymbol{\alpha}\I\, \mathbf{c} = \sfrac{1}{24}, 
&& \mathbf{b}\E\,\!^T\,\boldsymbol{\alpha}\E\,\boldsymbol{\gamma}\I\, \mathbf{c} = 0, \\
\label{eqn:GARK-ROW-IMEX-case1-order4}
& \mathbf{b}\E\,\!^T\,\boldsymbol{\alpha}\E\,\boldsymbol{\alpha}\I\, \mathbf{g}\I = 0, 
&& \mathbf{b}\E\,\!^T\,\boldsymbol{\alpha}\E\,\boldsymbol{\gamma}\I\, \mathbf{g}\I = 0, \\
& \mathbf{b}\I\,\!^T\,\boldsymbol{\alpha}\I\,\boldsymbol{\alpha}\E\, \mathbf{c} = \sfrac{1}{24}, 
&& \mathbf{b}\I\,\!^T\,\boldsymbol{\gamma}\I\,\boldsymbol{\alpha}\E\, \mathbf{c} = 0, \\
& \mathbf{b}\I\,\!^T\,\boldsymbol{\alpha}\I\,\boldsymbol{\alpha}\E\, \mathbf{g}\I = 0, 
&& \mathbf{b}\I\,\!^T\,\boldsymbol{\gamma}\I\,\boldsymbol{\alpha}\E\, \mathbf{g}\I = 0.
\end{alignedat}
\end{equation}
For $\Lb\I = \Jac_{n}\I$ the implicit part should be a Rosenbrock method of the desired order, the third order coupling conditions read:
\begin{equation}
\label{eqn:GARK-ROS-IMEX-case1-order3}
\mathbf{b}\E\,\!^T \, \boldsymbol{\alpha}\E\,\mathbf{g}\I = 0, \qquad
\mathbf{b}\I\,\!^T \, \boldsymbol{\beta}\I\,\mathbf{g}\I = 0,
\end{equation}
and the fourth coupling conditions are: 
\begin{equation}
\label{eqn:GARK-ROS-IMEX-case1-order4}
\begin{alignedat}{3}
&\mathbf{b}\E\,\!^T \, \big((\boldsymbol{\alpha}\E\,\mathbf{g}\I) \times \mathbf{c} \big) = 0, 
& \qquad &\mathbf{b}\E\,\!^T\,\boldsymbol{\alpha}\E\,\boldsymbol{\alpha}\E\, \mathbf{g}\I = 0, \\
&\mathbf{b}\E\,\!^T\,\boldsymbol{\alpha}\E\,\boldsymbol{\beta}\I\, \mathbf{c} = \sfrac{1}{24}, 
&&\mathbf{b}\E\,\!^T\,\boldsymbol{\alpha}\E\,\boldsymbol{\beta}\I\, \mathbf{g}\I = 0, \\
&\mathbf{b}\I\,\!^T\,\boldsymbol{\beta}\I\,\boldsymbol{\alpha}\E\, \mathbf{c} = \sfrac{1}{24}, 
&&\mathbf{b}\I\,\!^T\,\boldsymbol{\beta}\I\,\boldsymbol{\alpha}\E\, \mathbf{g}\I = 0.
\end{alignedat}
\end{equation}
\end{remark}

\begin{remark}
\label{rem:imex-same-b}
Another interesting special situation is when $\mathbf{b}\E=\mathbf{b}\I=\mathbf{b}$ in \cref{eqn:GARK-ROW-IMEX-case1}, in which case the scheme uses a single set of stages $\mathbf{k} = \mathbf{k}\I + \mathbf{k}\E$.
\end{remark}


The stability function \cref{eqn:GARK-ROS-stability} for an IMEX method \cref{eqn:GARK-ROW-IMEX-case1} becomes:

\begin{equation}
R= 1 + \begin{bmatrix} \b\E\,\!^T ~ \b\I\,\!^T \end{bmatrix} \, 
 \begin{bmatrix} z\E\,\!^{-1}\,\Idstage- \boldsymbol{\alpha}\EE & -\boldsymbol{\alpha}\EI \\ -\boldsymbol{\beta}\IE & z\I\,\!^{-1}\,\Idstage- \boldsymbol{\beta}\II \end{bmatrix}^{-1}\,\begin{bmatrix} \one_{s} \\ \one_{s} \end{bmatrix},
\end{equation}
where $s = s\E = s\I$. In the limit of infinite stiffness $z\I \to -\infty$:
\begin{equation*}
\begin{split}
&R
= R\I(\infty)+ z\E\,\left( \b\E\,\!^T - \b\I\,\!^T \boldsymbol{\beta}\II\,\!^{-1}\,\boldsymbol{\beta}\IE \right)
\mathbf{S}^{-1}\,\left( \Id - \boldsymbol{\alpha}\EI \boldsymbol{\beta}\II\,\!^{-1}\right)\,\one_{s}, \\
&\mathbf{S} = \Idstage- z\E\,(\boldsymbol{\alpha}\EE - \boldsymbol{\alpha}\EI\,\boldsymbol{\beta}\II\,\!^{-1}\,\boldsymbol{\beta}\IE).
 \end{split}
\end{equation*}
The second term is zero for stiffly accurate methods. Also this favorable situation arises when $\b\E = \b\I$ and $\boldsymbol{\beta}\IE = \boldsymbol{\beta}\II$.

\subsection{Implicit/linearly implicit GARK schemes}

%
The GARK-ROS framework allows to construct methods that are fully implicit in some partitions, and linearly implicit in other. For example, the explicit stage \cref{eqn:GARK-ROW-IMEX-stage} can be replaced by the following diagonally implicit stage (note the upper bound of the $\alpha_{i,j}\EE$ summation):
\begin{equation}
\label{eqn:GARK-ROW-IMLIM-scalar-stages-IM}
k_i\E = h\,\fun\E \mleft(\y_{n} +  \sum_{j=1}^{i} \alpha_{i,j}\EE\,k_j\E
+  \sum_{j=1}^{i-1 } \alpha_{i,j}\EI\,k_j\I \mright).
\end{equation}
The order conditions discussed above for the overall scheme remain unmodified.

\begin{remark}
By extension, one can construct GARK schemes that employ any combination of explicit, diagonally implicit, and linearly implicit methods to compute the stages associated with individual components.
\end{remark}

Moreover, one can formulate the stages \cref{eqn:GARK-ROW-IMLIM-scalar-stages-IM} as follows:
\begin{equation*}
\begin{split}
k_i\E &= h\,\fun\E \mleft(\y_{n} +  \sum_{j=1}^{i} \alpha_{i,j}\EE\,k_j\E
+  \sum_{j=1}^{i-1 } \alpha_{i,j}\EI\,k_j\I \mright), \\
&\quad + h \, \Lb\E\,\left( \sum_{j=1}^{i-1} \gamma_{i,j}\EE\,k_j\E + \sum_{j=1}^{i-1} \gamma_{i,j}\EI\,k_j\I \right). 
\end{split}
\end{equation*}
The computation remains explicit in $k_i\E$ when $\alpha_{i,i}\EE = 0$, and diagonally implicit when $\alpha_{i,i}\EE > 0$. The scheme no longer corresponds to either an explicit, or a diagonally implicit, GARK method. However, this formulation shows the power of the GARK-ROS framework to construct multimethods.

\section{Solution of index-1 differential-algebraic systems}
\label{sec:DAE-index-1}
%
Consider the singular perturbation problem \citep{Hairer_DAE_book,Hairer_book_II,Tikhonov_1985_ODE}
\begin{equation} 
\label{eqn:SPP}
\xx' = \fun(\xx,\zz), \qquad
\zz' = \varepsilon^{-1}\,\gun(\xx,\zz), 
\end{equation}
where $\varepsilon \ll 1$. The Jacobian $\gun_\zz$ is assumed to be invertible and with a negative logarithmic norm 
$\mu\left( \gun_\zz(\xx,\zz) \right) \le -1$ in an $\varepsilon$-independent neighborhood of the solution. 
Consequently, in the limit  $\varepsilon \to 0$ the system \cref{eqn:SPP} becomes an index-1 DAE \citep{Hairer_DAE_book,Hairer_book_II,Tikhonov_1985_ODE}:
\begin{equation}
\label{eqn:DAE-index1}
\xx' = \fun(\xx,\zz), \qquad
0 = \gun(\xx,\zz).
\end{equation}
The initial values $[\xx_n,\zz_n]$ are consistent if $\gun(\xx_n,\zz_n)=0$.
By the implicit function theorem the algebraic equation can be locally solved uniquely to obtain $\zz= \mathcal{G}(\xx)$.
Replacing this in the differential equation \cref{eqn:DAE-index1} leads to the following reduced ODE:
\begin{equation}
\label{eqn:reduced-ode}
\xx' = \fun\big(\xx,\mathcal{G}(\xx)\big) \eqqcolon \fun^{\textsc{red}}\big(\xx\big).
\end{equation}

Applying the GARK ROS scheme \cref{eqn:GARK-ROW-component} to \cref{eqn:SPP} gives:
\begin{subequations}
\label{eqn:GARK-ROW-component-SPP}
\begin{align}
\label{eqn:GARK-ROW-stage-componentX-SPP}
\mathbf{k} &= h\,\fun\Bigl(  \xx_{n}  + \boldsymbol{\alpha}\DiffDiff  \, \mathbf{k}, \zz_{n}  + \boldsymbol{\alpha}\DiffAlg  \, \boldsymbol{\ell} \Bigr)    
 + \nonumber \\
 & \quad + h\, \fun_{\xx | 0} \, \boldsymbol{\gamma}\DiffDiff \,  \mathbf{k}
+ h\,\fun_{\zz | 0} \,  \boldsymbol{\gamma}\DiffAlg \,  \boldsymbol{\ell},  \\ 
\label{eqn:GARK-ROW-stage-componentZ-SPP}
\boldsymbol{\ell} &= h\,\varepsilon^{-1}\, \gun\Bigl(  \xx_{n}  + \boldsymbol{\alpha}\AlgDiff  \, \mathbf{k}, \zz_{n}  + \boldsymbol{\alpha}\AlgAlg  \, \boldsymbol{\ell} \Bigr)   \\
\nonumber
& \quad + h\,\varepsilon^{-1}\, \gun_{\xx | 0} \, \boldsymbol{\gamma}\AlgDiff \,  \mathbf{k}
+ h\,\varepsilon^{-1}\,\gun_{\zz | 0} \,  \boldsymbol{\gamma}\AlgAlg \,  \boldsymbol{\ell},  \\
\label{eqn:GARK-ROW-solution-componentX-SPP}
\xx_{n+1} &= \xx_{n} + \b\Diff\,\!^T \, \mathbf{k}, \\
\label{eqn:GARK-ROW-solution-componentZ-SPP}
\zz_{n+1} &= \zz_{n} +  \b\Alg\,\!^T  \, \boldsymbol{\ell}.
\end{align}
\end{subequations} 
where, with a slight abuse of notation, we omit the explicit representation of the Kronecker products. The zero subscript means that the Jacobians are evaluated at the current step solution, e.g., $\gun_{\zz | 0} = \gun_{\zz}(\xx_{n},\zz_{n})$.

Taking the limit $\varepsilon \to 0$ changes \cref{eqn:GARK-ROW-stage-componentZ-SPP} into:
\label{eqn:GARK-ROW-component-DAE}
\begin{equation}
\label{eqn:GARK-ROW-stage-componentZ-DAE}
0 = \gun\Bigl(  \xx_{n}  + \boldsymbol{\alpha}\AlgDiff  \, \mathbf{k}, \zz_{n}  + \boldsymbol{\alpha}\AlgAlg  \, \boldsymbol{\ell} \Bigr) + \gun_{\xx | 0} \, \boldsymbol{\gamma}\AlgDiff \,  \mathbf{k} + \gun_{\zz | 0} \,  \boldsymbol{\gamma}\AlgAlg \,  \boldsymbol{\ell}.
\end{equation}
The $q$-th derivative of~\cref{eqn:GARK-ROW-stage-componentX-SPP} at $h=0$ is:
\begin{equation}
\label{eqn:q-derivative-k}
\begin{split}
\mathbf{k}^{(0)} & = 0; \\
 \mathbf{k}^{(1)} & = \fun(\xx_{n},\zz_{n}); \mbox{ and} \\
\mathbf{k}^{(q)} &= q\,\sum_{m+n \ge 2} \frac{\partial^{m+n}\fun}{\partial\xx^m \partial\zz^n} \Big|_{0}
\Big( \cdots , \boldsymbol{\alpha}\DiffDiff \, \mathbf{k}^{(\mu_i)} ,\cdots, \boldsymbol{\alpha}\DiffAlg \, \boldsymbol{\ell}^{(\nu_j)}, \cdots \Big)  \\
& \quad + q\, \fun_{\xx | 0} \, \boldsymbol{\beta}\DiffDiff \,  \mathbf{k}^{(q-1)} + q\,\fun_{\zz | 0} \,  \boldsymbol{\beta}\DiffAlg \,  \boldsymbol{\ell}^{(q-1)},  \\ 
& \quad \sum_{i=1}^m \mu_i + \sum_{j=1}^n \nu_i = q-1, \quad \textnormal{for }q \ge 2.
\end{split}
\end{equation}
Taking the $q$-th derivative of \cref{eqn:GARK-ROW-stage-componentZ-DAE} at $h=0$ gives:
\begin{equation}
\label{eqn:q-derivative-l}
\begin{split}
\nonumber
  0 &  =     \gun(\xx_{n},\zz_{n}); \\
\nonumber
  0 &  =     \boldsymbol{\beta}\AlgDiff \, \gun_{\xx | 0} \, \mathbf{k}^{(1)} +  \boldsymbol{\beta}\AlgAlg \,\gun_{\zz | 0} \,  \boldsymbol{\ell}^{(1)}; \mbox{ and} \\
0 &= \sum_{m+n \ge 2} \frac{\partial^{m+n}\gun}{\partial\xx^m \partial\zz^n} \Big|_{0}
\Big( \cdots , \boldsymbol{\alpha}\AlgDiff \, \mathbf{k}^{(\mu_i)} ,\cdots, \boldsymbol{\alpha}\AlgAlg \, \boldsymbol{\ell}^{(\nu_j)}, \cdots \Big)  \\
\nonumber
& \quad + \boldsymbol{\beta}\AlgDiff \, \gun_{\xx | 0} \, \mathbf{k}^{(q)} +  \boldsymbol{\beta}\AlgAlg \,\gun_{\zz | 0} \,  \boldsymbol{\ell}^{(q)}, \\
\nonumber
& \quad \sum_{i=1}^m \mu_i + \sum_{j=1}^n \nu_i = q, \quad \textnormal{for }q \ge 2.
\end{split}
\end{equation}
Using the notation $\boldsymbol{\omega}\AlgAlg = \boldsymbol{\beta}\AlgAlg\,\!^{-1}$ the second equation \cref{eqn:q-derivative-l} gives:
\begin{align*}
\boldsymbol{\ell}^{(q)} &= \boldsymbol{\omega}\AlgAlg \,(-\gun_{\zz | 0}^{-1}) \, \sum_{m+n \ge 2} \frac{\partial^{m+n}\gun}{\partial\xx^m \partial\zz^n} \Big|_{0}
\Big( \cdots , \boldsymbol{\alpha}\AlgDiff \, \mathbf{k}^{(\mu_i)} ,\cdots, \boldsymbol{\alpha}\AlgAlg \, \boldsymbol{\ell}^{(\nu_j)}, \cdots \Big),  \\
& \quad  + \boldsymbol{\omega}\AlgAlg \,\boldsymbol{\beta}\AlgDiff \, (-\gun_{\zz | 0}^{-1} \,\gun_{\xx | 0}) \, \mathbf{k}^{(q)}.
\end{align*}

We represent numerical solutions of GARK-ROW methods as NB-series over the set $\DAT$ of differential-algebraic trees \citep{Hairer_DAE_book,Hairer_book_II}. Let:
\begin{equation*}
\begin{split}
\mathbf{k} &= \textnormal{NB}\left( \boldsymbol{\theta}\Diff, [\xx_{n},\zz_{n}] \right), \quad
\boldsymbol{\ell} = \textnormal{NB}\left( \boldsymbol{\theta}\Alg, [\xx_{n},\zz_{n}] \right), \\
\xx_{n+1} &= \textnormal{NB}\left( \boldsymbol{\phi}\Diff, [\xx_{n},\zz_{n}] \right), \quad
\zz_{n+1} = \textnormal{NB}\left( \boldsymbol{\phi}\Alg, [\xx_{n},\zz_{n}] \right).
\end{split}
\end{equation*}
We have the following recurrences on NB-series coefficients:
\begin{equation*}
\begin{split}
\boldsymbol{\theta}\Diff(\utree) &= 0, \quad \forall\, \utree \in \DAT_\zz,  \\
\boldsymbol{\theta}\Diff(\tree) &= \begin{cases}
\boldsymbol{0}, & \tree = \emptyset, \\
\one, & \tree  = \tau_\xx,  \\
\Bigl( \Times_{i=1}^{m}  \boldsymbol{\alpha}\DiffDiff \boldsymbol{\theta}\Diff(\tree_i) \Bigr) \times
\Bigl( \Times_{j=1}^{n}  \boldsymbol{\alpha}\DiffAlg \boldsymbol{\theta}\Alg(\utree_j) \Bigr), \\
\qquad\qquad \tree=[\tree_1,\ldots,\tree_m,\utree_1,\ldots,\utree_n]_{\xx},
& m+n \ge 2, \\
\boldsymbol{\beta}\DiffDiff \boldsymbol{\theta}\Diff(\tree_1), & \tree=[\tree_1]_{\xx}, \\
\boldsymbol{\beta}\DiffAlg \boldsymbol{\theta}\Alg(\utree_1), & \tree=[\utree_1]_{\xx},
\end{cases}
\end{split}
\end{equation*}
and
\begin{equation*}
\begin{split}
\boldsymbol{\theta}\Alg(\tree) &= 0, \quad \forall\, \tree \in \DAT_\xx,  \\
\boldsymbol{\theta}\Alg(\utree) &= \begin{cases}
0, & \utree = \emptyset, \\
\boldsymbol{\omega}\AlgAlg\, \left( \big( \Times_{i=1}^{m}  \boldsymbol{\alpha}\AlgDiff \boldsymbol{\theta}\Diff(\tree_i) \big) \times
\big( \Times_{j=1}^{n}  \boldsymbol{\alpha}\AlgAlg \boldsymbol{\theta}\Alg(\utree_j) \big) \right), \\
\qquad\qquad\qquad \utree=[\tree_1,\ldots,\tree_m,\utree_1,\ldots,\utree_n]_{\zz}, & m+n \ge 2, \\
\boldsymbol{\omega}\AlgAlg \,\boldsymbol{\beta}\AlgDiff \, \boldsymbol{\theta}\Diff(\tree_1), & \utree=[\tree_1]_{\zz}.
\end{cases}
\end{split}
\end{equation*}

The final solutions \cref{eqn:GARK-ROW-solution-componentX-SPP} and  \cref{eqn:GARK-ROW-solution-componentZ-SPP} 
are represented, respectively, by NB-series with the following coefficients:
\begin{equation*}
\boldsymbol{\phi}\Diff(\tree) = \begin{cases}
1, & \tree = \emptyset, \\
\mathbf{b}\Diff\,\!^T\, \boldsymbol{\theta}\Diff(\tree),  & \textnormal{otherwise}.
\end{cases}
\qquad
\boldsymbol{\phi}\Alg(\utree) = \begin{cases}
1, & \utree = \emptyset, \\
\mathbf{b}\Alg\,\!^T\, \boldsymbol{\theta}\Alg(\utree),  & \textnormal{otherwise}.
\end{cases}
\end{equation*}
Equating the numerical and the exact solutions leads to the following.
\begin{theorem}[GARK-ROS order conditions for index-1 DAEs]
The numerical solution of the differential variable $\xx$ has order $p$ iff:
\begin{equation*}
\boldsymbol{\phi}\Diff(\tree) = \frac{1}{\gamma(\tree)} \quad \textnormal{for  } \tree \in \DAT_\xx,\quad \rho(\tree) \le p.
\end{equation*}
The numerical solution of the algebraic variable $\zz_n$ has order $q$ iff:
\begin{equation*}
\boldsymbol{\phi}\Alg(\utree) = \frac{1}{\gamma(\utree)} \quad \textnormal{for  } \utree \in \DAT_\zz,\quad \rho(\utree) \le q.
\end{equation*}
\end{theorem}

We form the stiff order conditions as follows: 
\begin{enumerate}
\item Meagre roots are labelled by $\mathbf{b}\Diff\,\!^T$ and fat roots by $\mathbf{b}\Alg\,\!^T\boldsymbol{\omega}\AlgAlg$.
\item A meagre node with a meagre parent is labelled $\boldsymbol{\alpha}\DiffDiff$ if it has multiple siblings, and by $\boldsymbol{\beta}\DiffDiff$ if it is the only child.
\item A meagre node with a fat parent is labelled $\boldsymbol{\alpha}\AlgDiff$ if it has multiple siblings, and by $\boldsymbol{\beta}\AlgDiff$ if it is the only child.
\item A fat node with a meagre parent is labelled $\boldsymbol{\alpha}\DiffAlg\,\boldsymbol{\omega}\AlgAlg$   if   it has multiple siblings, and $\boldsymbol{\beta}\DiffAlg\,\boldsymbol{\omega}\AlgAlg$ if it is the only child.
\item A fat node with a fat parent is labelled $\boldsymbol{\alpha}\AlgAlg\,\boldsymbol{\omega}\AlgAlg$ since it has multiple siblings.
\end{enumerate}
Based on this labelling, we form the stiff order conditions starting from the leaves and working toward the root:
\begin{enumerate}
\item Multiply the label of each leaf by $\one$ (of appropriate dimension).
\item Each node takes the component-wise product of its children's coefficients, and multiplies it by its label.
\end{enumerate}

\begin{remark}[Simplifying assumptions]
We make the simplifying assumption:
\begin{subequations}
\label{eqn:stiff-simplifying-assumptions}
\begin{equation}
\label{eqn:stiff-simplifying-assumption-a}
\boldsymbol{\beta}\AlgDiff = \boldsymbol{\beta}\AlgAlg \qquad \Rightarrow \qquad  \boldsymbol{\omega}\AlgAlg\,\boldsymbol{\beta}\AlgDiff = \Idstage.
\end{equation}
This assumption allows to simplify the order conditions as in \citep[Lemma 4.9, Section VI.4]{Hairer_book_II}. Order conditions for trees where a fat vertex is singly branched (by the structure of $\DAT$ trees, the child has to be meagre) involves products $\boldsymbol{\omega}\AlgAlg\, \boldsymbol{\beta}\AlgDiff$. The order conditions for such trees are redundant.
For example, \cref{eqn:stiff-simplifying-assumption-a} can be imposed when the scheme computes each $k_i\Diff$ before $k_i\Alg$. In this case one can have $\boldsymbol{\alpha}\AlgDiff$ and $\boldsymbol{\gamma}\AlgDiff$ lower triangular  (with non-zero diagonals), such that their sum matches $\boldsymbol{\beta}\AlgAlg$.

Note that when a singly branched meagre vertex is followed by a fat vertex  we have products $\boldsymbol{\beta}\DiffAlg \, \boldsymbol{\omega}\AlgAlg$. These trees are redundant when the following simplifying assumption holds:
\begin{equation}
\label{eqn:stiff-simplifying-assumption-b}
\boldsymbol{\beta}\DiffAlg = \boldsymbol{\beta}\AlgAlg \qquad \Rightarrow \qquad  \boldsymbol{\beta}\DiffAlg \, \boldsymbol{\omega}\AlgAlg = \Idstage.
\end{equation}
\end{subequations}
For example, \cref{eqn:stiff-simplifying-assumption-b} can be imposed when the scheme computes each $k_i\Alg$ before $k_i\Diff$. In this case one can have $\boldsymbol{\alpha}\DiffAlg$ and $\boldsymbol{\gamma}\DiffAlg$ lower triangular (with non-zero diagonals), such that their sum matches $\boldsymbol{\beta}\AlgAlg$.

However, imposing both conditions \cref{eqn:stiff-simplifying-assumption-a} and \cref{eqn:stiff-simplifying-assumption-b} leads to the requirement that  $k_i\Alg$ and $k_i\Diff$ are computed together, therefore the resulting scheme is no longer decoupled. Stiff order conditions for Rosenbrock methods, which compute a single set of stages, benefit from both conditions \cref{eqn:stiff-simplifying-assumptions} \citep{Hairer_book_II}.
\end{remark}

2

\begin{table}
\centering
\begin{tabular}{|| C{1.75em} || C{13em} || C{14em}  ||  C{1.75em} ||}
\hline\hline
$\tree$ &  Labels & $\phi(\tree)$ & $\gamma(\tree)$ \\
\hline\hline
$\utree_{2,1}$ &
\begin{tikzpicture}[grow=up, level distance=0.6cm]\Tree [.\node[circle,draw,thick,opacity=1,fill=white,inner sep=0,minimum size=4,label=right:$\mathbf{b}\Alg\,\!^T \boldsymbol{\omega}\AlgAlg$](j){}; [.\node[circle,draw,thick,opacity=1,fill=black,inner sep=0,minimum size=4,label=right:$\boldsymbol{\alpha}\AlgDiff$](l){};  ][.\node[circle,draw,thick,opacity=1,fill=black,inner sep=0,minimum size=4,label=right:$\boldsymbol{\alpha}\AlgDiff$](k){};  ] ]\end{tikzpicture} &   
$\mathbf{b}\Alg\,\!^T \boldsymbol{\omega}\AlgAlg \, \mathbf{c}\AlgDiff\,\!^{\times 2}$ & 1 \\ 
\hline
$\utree_{3,1}$ &
\begin{tikzpicture}[grow=up, level distance=0.6cm]\Tree [.\node[circle,draw,thick,opacity=1,fill=white,inner sep=0,minimum size=4,label=right:$\mathbf{b}\Alg\,\!^T \boldsymbol{\omega}\AlgAlg$](j){}; [.\node[circle,draw,thick,opacity=1,fill=black,inner sep=0,minimum size=4,label=right:$\boldsymbol{\alpha}\AlgDiff$](m){};  ][.\node[circle,draw,thick,opacity=1,fill=black,inner sep=0,minimum size=4,label=right:$\boldsymbol{\alpha}\AlgDiff$](l){};  ][.\node[circle,draw,thick,opacity=1,fill=black,inner sep=0,minimum size=4,label=right:$\boldsymbol{\alpha}\AlgDiff$](k){};  ] ]\end{tikzpicture} & 
$\mathbf{b}\Alg\,\!^T \boldsymbol{\omega}\AlgAlg \, \boldsymbol{c}\AlgDiff\,\!^{\times 3}$ & 1 \\ 
\hline
$\utree_{3,2}$ &
\begin{tikzpicture}[grow=up, level distance=0.6cm]\Tree [.\node[circle,draw,thick,opacity=1,fill=white,inner sep=0,minimum size=4,label=right:$\mathbf{b}\Alg\,\!^T \boldsymbol{\omega}\AlgAlg$](j){}; [.\node[circle,draw,thick,opacity=1,fill=black,inner sep=0,minimum size=4,label=right:$\boldsymbol{\alpha}\AlgDiff$](l){};  ][.\node[circle,draw,thick,opacity=1,fill=black,inner sep=0,minimum size=4,label=right:$\boldsymbol{\alpha}\AlgDiff$](k){}; [.\node[circle,draw,thick,opacity=1,fill=black,inner sep=0,minimum size=4,label=right:$\boldsymbol{\beta}\DiffDiff$](m){};  ] ] ]\end{tikzpicture} &  
$\mathbf{b}\Alg\,\!^T \boldsymbol{\omega}\AlgAlg \, \big( (\boldsymbol{\alpha}\AlgDiff\,\mathbf{e}\DiffDiff) \times \mathbf{c}\AlgDiff\big)$ & 2 \\ 
\hline
$\utree_{3,3}$ & \begin{tikzpicture}[grow=up, level distance=0.6cm]\Tree [.\node[circle,draw,thick,opacity=1,fill=white,inner sep=0,minimum size=4,label=right:$\mathbf{b}\Alg\,\!^T \boldsymbol{\omega}\AlgAlg$](j){}; [.\node[circle,draw,thick,opacity=1,fill=black,inner sep=0,minimum size=4,label=right:$\boldsymbol{\alpha}\AlgDiff$](l){};  ][.\node[circle,draw,thick,opacity=1,fill=white,inner sep=0,minimum size=4,label=right:$\boldsymbol{\alpha}\AlgAlg\, \boldsymbol{\omega}\AlgAlg$](k){}; [.\node[circle,draw,thick,opacity=1,fill=black,inner sep=0,minimum size=4,label=right:$\boldsymbol{\alpha}\AlgDiff$](n){};  ][.\node[circle,draw,thick,opacity=1,fill=black,inner sep=0,minimum size=4,label=right:$\boldsymbol{\alpha}\AlgDiff$](m){};  ] ] ]\end{tikzpicture}  & 
$\begin{array}{c} \mathbf{b}\Alg\,\!^T \boldsymbol{\omega}\AlgAlg \cdot \\ \cdot \Big( \big( \boldsymbol{\alpha}\AlgAlg\, \boldsymbol{\omega}\AlgAlg\, \mathbf{c}\AlgDiff\,\!^{\times 2} \big) \\ \times \mathbf{c}\AlgDiff\Big) \end{array}$ & 1 \\  
\hline\hline
$\tree_{3,1}$ &
\begin{tikzpicture}[grow=up, level distance=0.6cm]\Tree [.\node[circle,draw,thick,opacity=1,fill=black,inner sep=0,minimum size=4,label=right:$\mathbf{b}\Diff\,\!^T$](j){}; [.\node[circle,draw,thick,opacity=1,fill=white,inner sep=0,minimum size=4,label=right:$\boldsymbol{\beta}\DiffAlg\,\boldsymbol{\omega}\AlgAlg$](k){}; [.\node[circle,draw,thick,opacity=1,fill=black,inner sep=0,minimum size=4,label=right:$\boldsymbol{\alpha}\AlgDiff$](m){};  ][.\node[circle,draw,thick,opacity=1,fill=black,inner sep=0,minimum size=4,label=right:$\boldsymbol{\alpha}\AlgDiff$](l){};  ] ] ]\end{tikzpicture}  & 
$\mathbf{b}\Diff\,\!^T \, \boldsymbol{\beta}\DiffAlg\,\boldsymbol{\omega}\AlgAlg\, \mathbf{c}\AlgDiff\,\!^{\times 2}$ & 3 \\ 
\hline
$\tree_{4,1}$ & \begin{tikzpicture}[grow=up, level distance=0.6cm]\Tree [.\node[circle,draw,thick,opacity=1,fill=black,inner sep=0,minimum size=4,label=right:$\mathbf{b}\Diff\,\!^T$](j){}; [.\node[circle,draw,thick,opacity=1,fill=black,inner sep=0,minimum size=4,label=right:$\boldsymbol{\alpha}\DiffDiff$](l){};  ][.\node[circle,draw,thick,opacity=1,fill=white,inner sep=0,minimum size=4,label=right:$\boldsymbol{\alpha}\DiffAlg\, \boldsymbol{\omega}\AlgAlg$](k){}; [.\node[circle,draw,thick,opacity=1,fill=black,inner sep=0,minimum size=4,label=right:$\boldsymbol{\alpha}\AlgDiff$](n){};  ][.\node[circle,draw,thick,opacity=1,fill=black,inner sep=0,minimum size=4,label=right:$\boldsymbol{\alpha}\AlgDiff$](m){};  ] ] ]\end{tikzpicture}  & 
$\begin{array}{c} \mathbf{b}\Diff\,\!^T \,\Big(  \mathbf{c}\DiffDiff \times \\ \big( \boldsymbol{\alpha}\DiffAlg\,  \boldsymbol{\omega}\AlgAlg\,\mathbf{c}\AlgDiff\,\!^{\times 2} \big) \Big)\end{array}$ & 4 \\  
\hline
$\tree_{4,2}$ & \begin{tikzpicture}[grow=up, level distance=0.6cm]\Tree [.\node[circle,draw,thick,fill=black,inner sep=0,minimum size=4,opacity=1,label=right:$\mathbf{b}\Diff\,\!^T$](j){}; [.\node[circle,draw,thick,opacity=1,fill=white,inner sep=0,minimum size=4,label=right:$\boldsymbol{\beta}\DiffAlg\,\boldsymbol{\omega}\AlgAlg$](k){}; [.\node[circle,draw,thick,opacity=1,fill=black,inner sep=0,minimum size=4,label=right:$\boldsymbol{\alpha}\AlgDiff$](n){};  ][.\node[circle,draw,thick,opacity=1,fill=black,inner sep=0,minimum size=4,label=right:$\boldsymbol{\alpha}\AlgDiff$](m){};  ][.\node[circle,draw,thick,opacity=1,fill=black,inner sep=0,minimum size=4,label=right:$\boldsymbol{\alpha}\AlgDiff$](l){};  ] ] ]\end{tikzpicture}  & 
$\mathbf{b}\Diff\,\!^T \, \boldsymbol{\beta}\DiffAlg\,\boldsymbol{\omega}\AlgAlg\, \mathbf{c}\AlgDiff\,\!^{\times 3}$ & 4 \\ 
\hline
$\tree_{4,3}$ & \begin{tikzpicture}[grow=up, level distance=0.6cm]\Tree [.\node[circle,draw,thick,opacity=1,fill=black,inner sep=0,minimum size=4,label=right:$\mathbf{b}\Diff\,\!^T$](j){}; [.\node[circle,draw,thick,opacity=1,fill=white,inner sep=0,minimum size=4,label=right:$\boldsymbol{\beta}\DiffAlg\boldsymbol{\omega}\AlgAlg$](k){}; [.\node[circle,draw,thick,opacity=1,fill=black,inner sep=0,minimum size=4,label=right:$\boldsymbol{\alpha}\AlgDiff$](m){};  ][.\node[circle,draw,thick,opacity=1,fill=black,inner sep=0,minimum size=4,label=right:$\boldsymbol{\alpha}\AlgDiff$](l){}; [.\node[circle,draw,thick,opacity=1,fill=black,inner sep=0,minimum size=4,label=right:$\boldsymbol{\beta}\DiffDiff$](n){};  ] ] ] ]\end{tikzpicture}  & 
$\mathbf{b}\Diff\,\!^T \,\boldsymbol{\beta}\DiffAlg\,\boldsymbol{\omega}\AlgAlg\, \bigl( \mathbf{c}\AlgDiff \times (\boldsymbol{\alpha}\AlgDiff\,\mathbf{e}\DiffDiff) \bigr) $ & 8 \\ 
\hline
$\tree_{4,4}$ & \begin{tikzpicture}[grow=up, level distance=0.6cm]\Tree [.\node[circle,draw,thick,opacity=1,fill=black,inner sep=0,minimum size=4,label=right:$\mathbf{b}\Diff\,\!^T$](j){}; [.\node[circle,draw,thick,opacity=1,fill=black,inner sep=0,minimum size=4,label=right:$\boldsymbol{\beta}\DiffDiff$](k){}; [.\node[circle,draw,thick,opacity=1,fill=white,inner sep=0,minimum size=4,label=right:$\boldsymbol{\beta}\DiffAlg\boldsymbol{\omega}\AlgAlg$](l){}; [.\node[circle,draw,thick,opacity=1,fill=black,inner sep=0,minimum size=4,label=right:$\boldsymbol{\alpha}\AlgDiff$](n){};  ][.\node[circle,draw,thick,opacity=1,fill=black,inner sep=0,minimum size=4,label=right:$\boldsymbol{\alpha}\AlgDiff$](m){};  ] ] ] ]\end{tikzpicture}  & 
$\begin{array}{c} \mathbf{b}\Diff\,\!^T \,\boldsymbol{\beta}\DiffDiff\,\boldsymbol{\beta}\DiffAlg \cdot \\ \cdot \boldsymbol{\omega}\AlgAlg \mathbf{c}\AlgDiff\,\!^{\times 2} \end{array}$ & 12 \\ 
\hline\hline
\end{tabular}
\caption{$\DAT$ trees and order conditions for GARK-ROS numerical solution using the simplifying assumption \cref{eqn:stiff-simplifying-assumption-a}. Follows \citep[Table 4.1, Section VI.4]{Hairer_book_II}.}
\label{tab:Ros-DAT-trees_numeric}
\end{table}



Following \citep[Table 4.1, Section VI.4]{Hairer_book_II}, the first $\DAT$ trees are shown  in Table  \ref{tab:Ros-DAT-trees_numeric}.   Only the trees remaining after the simplifying assumption \cref{eqn:stiff-simplifying-assumption-a} is imposed are shown. We have the following result.

\begin{theorem}[Algebraic order conditions for index-1 DAE solution]
The algebraic order conditions are as follows. 
\begin{subequations}
\label{eqn:index-1-order}
\begin{align}
\label{eqn:index-1-order-z2}
&\textnormal{order 2 }(\zz):
\begin{cases}
\mathbf{b}\Alg\,\!^T \boldsymbol{\omega}\AlgAlg \, \mathbf{c}\AlgDiff\,\!^{\times 2} = 1;
\end{cases}
\\
\label{eqn:index-1-order-z3}
&\textnormal{order 3 }(\zz):
\begin{cases}
\mathbf{b}\Alg\,\!^T \boldsymbol{\omega}\AlgAlg \, \boldsymbol{c}\AlgDiff\,\!^{\times 3} = 1, \\
\mathbf{b}\Alg\,\!^T \boldsymbol{\omega}\AlgAlg \, \big( (\boldsymbol{\alpha}\AlgDiff\,\boldsymbol{e}\DiffDiff) \times \mathbf{c}\AlgDiff\big) = \frac{1}{2}, \\
\mathbf{b}\Alg\,\!^T \boldsymbol{\omega}\AlgAlg \,\Big( \big( \boldsymbol{\alpha}\AlgAlg\, \boldsymbol{\omega}\AlgAlg\, \mathbf{c}\AlgDiff\,\!^{\times 2} \big) \times \mathbf{c}\AlgDiff\Big) =  1;
\end{cases}
\\
\label{eqn:index-1-order-x3}
&\textnormal{order 3 }(\xx):
\begin{cases}
\mathbf{b}\Diff\,\!^T \, \boldsymbol{\beta}\DiffAlg\,\boldsymbol{\omega}\AlgAlg\, \mathbf{c}\AlgDiff\,\!^{\times 2} = \frac{1}{3};
\end{cases}
\\
\label{eqn:index-1-order-x4}
&\textnormal{order 4 }(\xx):
\begin{cases}
\mathbf{b}\Diff\,\!^T \,\Big( \big( \boldsymbol{\alpha}\DiffAlg\,  \boldsymbol{\omega}\AlgAlg\,\mathbf{c}\AlgDiff\,\!^{\times 2} \big) \times \mathbf{c}\DiffDiff \Big) =  \frac{1}{4}, \\
\mathbf{b}\Diff\,\!^T \, \boldsymbol{\beta}\DiffAlg\,\boldsymbol{\omega}\AlgAlg\, \mathbf{c}\AlgDiff\,\!^{\times 3} =  \frac{1}{4}, \\
\mathbf{b}\Diff\,\!^T \,\boldsymbol{\beta}\DiffAlg\,\boldsymbol{\omega}\AlgAlg\, \bigl( \mathbf{c}\AlgDiff \times (\boldsymbol{\alpha}\AlgDiff\,\mathbf{e}\DiffDiff) \bigr) = \frac{1}{8}, \\
\mathbf{b}\Diff\,\!^T \,\boldsymbol{\beta}\DiffDiff\,\boldsymbol{\beta}\DiffAlg\boldsymbol{\omega}\AlgAlg \mathbf{c}\AlgDiff\,\!^{\times 2} = \frac{1}{12}.
\end{cases}
\end{align}
\end{subequations}
\end{theorem}

\begin{remark}[Special case IMEX method]
For the IMEX GARK scheme with the special structure discussed in Remark \ref{rem:imex-special-case} the order conditions are as follows. 
The algebraic order conditions for $\zz$ are the ones of the implicit component. Thus, if the implicit component has index-1 DAE order $q$ for $\zz$ then the IMEX GARK component inherits this property.
The index-1 DAE conditions for $\xx$ are, for order three:
\begin{equation} \label{eqn:GARK-ROS-IMEX-DAE-case1-order3}
\mathbf{b}\Diff\,\!^T \, \boldsymbol{\alpha}\Diff\,\boldsymbol{\omega}\Alg\, \mathbf{c}^{\times 2} = \sfrac{1}{3}, 
\end{equation}
and for order four:
\begin{equation} 
\label{eqn:GARK-ROS-IMEX-DAE-case1-order4}
\begin{split}
\mathbf{b}\Diff\,\!^T \,\Big( \big( \boldsymbol{\alpha}\Diff\,  \boldsymbol{\omega}\Alg\,\mathbf{c}^{\times 2} \big) \times \mathbf{c} \Big) &=  \sfrac{1}{4}, \qquad
\mathbf{b}\Diff\,\!^T \, \boldsymbol{\alpha}\Diff\,\boldsymbol{\omega}\Alg\, \mathbf{c}^{\times 3} =  \sfrac{1}{4}, \\
\mathbf{b}\Diff\,\!^T \,\boldsymbol{\alpha}\Diff\,\boldsymbol{\omega}\Alg\, \bigl( \mathbf{c} \times (\boldsymbol{\alpha}\Alg\,\mathbf{c}) \bigr) &=  \sfrac{1}{8}, \qquad
\mathbf{b}\Diff\,\!^T \,\boldsymbol{\alpha}\Diff\,\boldsymbol{\alpha}\Diff\boldsymbol{\omega}\Alg \mathbf{c}^{\times 2} =  \sfrac{1}{12}.
\end{split}
\end{equation}
They are solved together with the classical order conditions \cref{eqn:GARK-ROW-IMEX-case1-order3} and 
\cref{eqn:GARK-ROW-IMEX-case1-order4}.

\end{remark}

\begin{remark}[Order conditions for inconsistent initial values]
Incon- \linebreak sistent initial conditions $\gun(\xx_n,\zz_n) \ne 0$ lead to additional error terms in the numerical solution \citep[Table 4.2, Section VI.4]{Hairer_book_II}. These error terms correspond to solution derivatives that contain $-\gun_{\zz | 0}^{-1}\, \gun(\xx_n,\zz_n)$ terms, and therefore to $\DAT$ trees that have fat leaves. Assume that the inconsistency satisfies:
\[
\Vert -\gun_{\zz | 0}^{-1}\, \gun(\xx_n,\zz_n) \Vert \le \delta.
\]
\ifreport
The first $\DAT$ trees with fat leaves are summarized in \cref{tab:Ros-DAT-IC-trees_numeric}, which follows 
\citep[Table 4.2, Section VI.4]{Hairer_book_II}. 
\fi
Each tree corresponds to an error term due to the initial value inconsistency; the number of fat leaves gives the power of $\delta$ , and the number of meagre nodes the power of $h$ in the corresponding error term. 

\ifreport
The order conditions are given in \cref{tab:Ros-DAT-IC-trees_numeric}.
\fi
Let $\boldsymbol{o}\Alg \coloneqq \boldsymbol{\omega}\AlgAlg \, \one\Alg$. The first order conditions for $\zz$ read:
\begin{subequations}
\label{eqn:inconsistent-ic-z}
\begin{align}
\label{eqn:inconsistent-ic-z-1}
\mathcal{O}(\delta): & ~~ \mathbf{b}\Alg\,\!^T\,\boldsymbol{o}\Alg = 1,  \\ 
\label{eqn:inconsistent-ic-z-2}
\mathcal{O}(h\delta): & ~~ \mathbf{b}\Alg\,\!^T \boldsymbol{\omega}\AlgAlg \cdot \left(\mathbf{c}\AlgDiff \times \boldsymbol{\alpha}\AlgAlg\, \boldsymbol{o}\Alg \right)  = 1,
\end{align}
\end{subequations}
and the first ones for $\xx$ are:
\begin{subequations}
\label{eqn:inconsistent-ic-x}
\begin{align}
\mathcal{O}(h\delta): &~~ \mathbf{b}\Diff\,\!^T \, \boldsymbol{\beta}\DiffAlg\, \boldsymbol{o}\Alg = 1, \label{eqn:inconsistent-ic-x-1} \\ 
\mathcal{O}(h^2\delta): &~~ \mathbf{b}\Diff\,\!^T \, (\mathbf{c}\DiffDiff \times \boldsymbol{\alpha}\DiffAlg\,\boldsymbol{o}\Alg) = \sfrac{1}{2}, \\ 
\mathcal{O}(h^2\delta): &~~ \mathbf{b}\Diff\,\!^T \, \boldsymbol{\beta}\DiffDiff\,\boldsymbol{\beta}\DiffAlg\,\boldsymbol{o}\Alg = \sfrac{1}{2}, \\ 
\mathcal{O}(h^2\delta): &~~ \mathbf{b}\Diff\,\!^T \, \boldsymbol{\beta}\DiffAlg\,\boldsymbol{\omega}\AlgAlg\cdot \left( \mathbf{c}\AlgDiff \times \boldsymbol{\alpha}\AlgAlg \boldsymbol{o}\Alg \right) = \sfrac{1}{2}. 
\end{align}
\end{subequations}
If the numerical solution satisfies all the additional order conditions \cref{eqn:inconsistent-ic-z} and \cref{eqn:inconsistent-ic-x} then the (additional) local error in $\xx$ due to inconsistent initial conditions is $\mathcal{O}(h^3\delta + h \delta^2)$, and the local error in  $\zz$ is $\mathcal{O}(h^2\delta + \delta^2)$. 
\end{remark}

\ifreport

\begin{table}[ht!]
	\centering
	\begin{tabular}{|| C{2.5em} || C{11em} || C{14em}  ||  C{1.75em} ||}
		\hline\hline
		Error &  Labels & $\phi(\cdot)$ & $\gamma(\cdot)$ \\
		\hline\hline
		$\begin{array}{c} \delta \\ (\zz) \end{array}$ &
		\begin{tikzpicture}[grow=up, level distance=0.6cm]\Tree [.\node[circle,draw,thick,opacity=1,fill=white,inner sep=0,minimum size=4,label=right:$\mathbf{b}\Alg\,\!^T\,\boldsymbol{\omega}\AlgAlg$](b){};  ]\end{tikzpicture} 
		&   $\mathbf{b}\Alg\,\!^T\,\boldsymbol{\omega}\AlgAlg \, \one\Alg$ & 1 \\ 
		\hline
		$\begin{array}{c} h \delta \\ (\zz) \end{array}$ &
		\begin{tikzpicture}[grow=up, level distance=0.6cm]\Tree [.\node[circle,draw,thick,opacity=1,fill=white,inner sep=0,minimum size=4,label=right:$\mathbf{b}\Alg\,\!^T \boldsymbol{\omega}\AlgAlg$](j){}; [.\node[circle,draw,thick,opacity=1,fill=white,inner sep=0,minimum size=4,label=right:$\boldsymbol{\alpha}\AlgAlg \boldsymbol{\omega}\AlgAlg $](l){};  ][.\node[circle,draw,thick,opacity=1,fill=black,inner sep=0,minimum size=4,label=right:$\boldsymbol{\alpha}\AlgDiff$](k){};  ] ]\end{tikzpicture} &   
		$\begin{array}{c} \mathbf{b}\Alg\,\!^T \boldsymbol{\omega}\AlgAlg \cdot \\ \left(\mathbf{c}\AlgDiff \times \boldsymbol{\alpha}\AlgAlg \boldsymbol{\omega}\AlgAlg \one\Alg \right) \end{array}$ & 1 \\ 
		\hline\hline
		$\begin{array}{c} h \delta \\ (\xx) \end{array}$ &
		\begin{tikzpicture}[grow=up, level distance=0.6cm]\Tree [.\node[circle,draw,thick,opacity=1,fill=black,inner sep=0,minimum size=4,label=right:$\mathbf{b}\Diff\,\!^T$](j){}; [.\node[circle,draw,thick,opacity=1,fill=white,inner sep=0,minimum size=4,label=right:$\boldsymbol{\beta}\DiffAlg\,  \boldsymbol{\omega}\AlgAlg$](k){};  ] ]\end{tikzpicture} 
		&    $\mathbf{b}\Diff\,\!^T \, \boldsymbol{\beta}\DiffAlg\,  \boldsymbol{\omega}\AlgAlg \, \one\Alg$ & 1 \\ 
		\hline
		$\begin{array}{c} h^2 \delta \\ (\xx) \end{array}$ &
		\begin{tikzpicture}[grow=up, level distance=0.6cm]\Tree [.\node[circle,draw,thick,opacity=1,fill=black,inner sep=0,minimum size=4,label=right:$\mathbf{b}\Diff\,\!^T$](j){}; [.\node[circle,draw,thick,opacity=1,fill=white,inner sep=0,minimum size=4,label=right:$\boldsymbol{\alpha}\DiffAlg \boldsymbol{\omega}\AlgAlg$](l){};  ][.\node[circle,draw,thick,opacity=1,fill=black,inner sep=0,minimum size=4,label=right:$\boldsymbol{\alpha}\DiffDiff$](k){};  ] ]\end{tikzpicture} &   
		$\mathbf{b}\Diff\,\!^T \, (\mathbf{c}\DiffDiff \times \boldsymbol{\alpha}\DiffAlg\,\boldsymbol{\omega}\AlgAlg \, \one\Alg)$ & 2 \\ 
		\hline
		$\begin{array}{c} h^2 \delta \\ (\xx) \end{array}$ &
		\begin{tikzpicture}[grow=up, level distance=0.6cm]\Tree [.\node[circle,draw,thick,opacity=1,fill=black,inner sep=0,minimum size=4,label=right:$\mathbf{b}\Diff\,\!^T$](j){}; [.\node[circle,draw,thick,opacity=1,fill=black,inner sep=0,minimum size=4,label=right:$\boldsymbol{\beta}\DiffDiff$](k){}; [.\node[circle,draw,thick,opacity=1,fill=white,inner sep=0,minimum size=4,label=right:$\boldsymbol{\beta}\DiffAlg\,\boldsymbol{\omega}\AlgAlg$](l){};  ] ] ]\end{tikzpicture} &
		$\mathbf{b}\Diff\,\!^T \, \boldsymbol{\beta}\DiffDiff\,\boldsymbol{\beta}\DiffAlg\,\boldsymbol{\omega}\AlgAlg\,\one\Alg$ & 2 \\ 
		\hline
		%
		$\begin{array}{c} h^2 \delta \\ (\xx) \end{array}$ &
		\begin{tikzpicture}[grow=up, level distance=0.6cm]\Tree [.\node[circle,draw,thick,opacity=1,fill=black,inner sep=0,minimum size=4,label=right:$\mathbf{b}\Diff\,\!^T$](j){}; [.\node[circle,draw,thick,opacity=1,fill=white,inner sep=0,minimum size=4,label=right:$\boldsymbol{\beta}\DiffAlg\,\boldsymbol{\omega}\AlgAlg$](k){}; [.\node[circle,draw,thick,opacity=1,fill=white,inner sep=0,minimum size=4,label=right:$\boldsymbol{\alpha}\AlgAlg \boldsymbol{\omega}\AlgAlg$](m){};  ][.\node[circle,draw,thick,opacity=1,fill=black,inner sep=0,minimum size=4,label=right:$\boldsymbol{\alpha}\AlgDiff$](l){};  ] ] ]\end{tikzpicture}  & 
		$\begin{array}{l}\mathbf{b}\Diff\,\!^T \, \boldsymbol{\beta}\DiffAlg\,\boldsymbol{\omega}\AlgAlg\cdot \\\left( \mathbf{c}\AlgDiff \times \boldsymbol{\alpha}\AlgAlg \boldsymbol{\omega}\AlgAlg \one\Alg\right)\end{array}$ & 2 \\ 
		\hline\hline
	\end{tabular}
	\caption{Trees and order conditions for inconsistent initial conditions in GARK-ROS index-1 DAE numerical solution, using the simplifying assumption \cref{eqn:stiff-simplifying-assumption-a}. Follows \citep[Table 4.2, Section VI.4]{Hairer_book_II}.}
	\label{tab:Ros-DAT-IC-trees_numeric}
\end{table}

\fi



\section{Practical GARK-ROS/ROW methods}
\label{sec:methods}

In this section, we develop new linearly implicit GARK methods up to order four.  All of the IMEX-type methods are derived using the classical order conditions from \cref{sec:IMEX-GARK-ROW}, and the stiff order conditions from \cref{sec:DAE-index-1} where appropriate. We note that the coefficients we derive also satisfy the corresponding decoupled IMIM order conditions from \cref{sec:decoupled-GARK-ROW} with
\begin{equation*}
    \boldsymbol{\alpha} = \boldsymbol{\alpha}\I,
    \quad
    \underline{\boldsymbol{\alpha}} = \overline{\boldsymbol{\alpha}} = \boldsymbol{\alpha}\E,
    \quad
    \boldsymbol{\gamma} = \boldsymbol{\gamma}\I,
    \quad
    \underline{\boldsymbol{\gamma}} = \overline{\boldsymbol{\gamma}} = \boldsymbol{\gamma}\E.
\end{equation*}
Therefore, they can be used to construct decoupled IMIM GARK-ROW/GARK-ROS schemes of the same order.

\subsection{Second order implicit/linearly implicit/explicit multimethod} 

Consider the system \cref{eqn:additively-partitioned-ode} with $\nparts=3$ partitions where the first partition is nonstiff and the other two are stiff.  To showcase the flexibility of the linearly implicit GARK framework, we develop a second order multimethod that combines an explicit Runge--Kutta method, an implicit Runge--Kutta method, with a Rosenbrock method.  In particular, we use the implicit and explicit trapezoidal rules:
\begin{equation*}
    \begin{butchertableau}{c|c}
        \c^{\textsc{it}} & \A^{\textsc{it}} \\ \hline
        & (\b^{\textsc{it}})^T
    \end{butchertableau} = \raisebox{8pt}{$\begin{butchertableau}{c|cc}
        {\scriptstyle 0}  & {\scriptstyle 0} & {\scriptstyle 0} \\
        {\scriptstyle 1}  & \frac{1}{2} & \frac{1}{2} \\ \hline
        & \frac{1}{2} & \frac{1}{2}
    \end{butchertableau}$},
    \qquad
    \begin{butchertableau}{c|c}
        \c^{\textsc{et}} & \A^{\textsc{et}} \\ \hline
        & (\b^{\textsc{et}})^T
    \end{butchertableau} = \raisebox{8pt}{$\begin{butchertableau}{c|cc}
        {\scriptstyle 0}  & {\scriptstyle 0} & {\scriptstyle 0} \\
        {\scriptstyle 1} & {\scriptstyle 1} & {\scriptstyle 0} \\ \hline
        & \frac{1}{2} & \frac{1}{2}
    \end{butchertableau}$},
\end{equation*}
as well as the stiffly accurate, L-stable Rosenbrock scheme with coefficients
\begin{equation*}
    \boldsymbol{\alpha}^{\textsc{ros2}} = \begin{bmatrix}
        {\scriptstyle 0} & {\scriptstyle 0} \\
        {\scriptstyle 1} & {\scriptstyle 0} \\
    \end{bmatrix},
    \quad
    \boldsymbol{\gamma}^{\textsc{ros2}} = \begin{bmatrix}
        \hphantom{-}\gamma & {\scriptstyle 0} \\
         -\gamma & \gamma
    \end{bmatrix},
    \quad
    \b^{\textsc{ros2}} = \begin{bmatrix}
        1 - \gamma & \gamma
    \end{bmatrix}^T, 
    \quad
    \gamma = 1 - \frac{\sqrt{2}}{2}.
\end{equation*}

There are six $\boldsymbol{\alpha}$ coupling matrices and two $\boldsymbol{\gamma}$ coupling matrices to be determined for this multimethod, which offers numerous degrees of freedom.  We use the simplifying assumptions of \cref{rem:imex-special-case} with a slight modification to ensure the fully implicit and linearly implicit stages are decoupled.  The linearly implicit GARK scheme defined by the tableau
\begin{equation*}
    \begin{butchertableau}{c|c}
        \A & \G \\ \hline
        \b^T &
    \end{butchertableau}
    ~=~
    ~\raisebox{15.5pt}{$\begin{butchertableau}{ccc|ccc}
        \A^{\textsc{et}} & \A^{\textsc{et}} & \A^{\textsc{et}} & \mathbf{0} & \mathbf{0} & \mathbf{0} \\
        \A^{\textsc{it}} & \A^{\textsc{it}} & \A^{\textsc{et}} & \mathbf{0} & \mathbf{0} & \mathbf{0} \\
        \boldsymbol{\alpha}^{\textsc{ros2}} & \boldsymbol{\alpha}^{\textsc{ros2}} & \boldsymbol{\alpha}^{\textsc{ros2}} & \boldsymbol{\gamma}^{\textsc{ros2}} & \boldsymbol{\gamma}^{\textsc{ros2}} & \boldsymbol{\gamma}^{\textsc{ros2}} \\ \hline
        (\b^{\textsc{et}})^T & (\b^{\textsc{it}})^T & (\b^{\textsc{ros2}})^T
    \end{butchertableau}$}
\end{equation*}
maintains the second order of the base methods and is suitable for index 1 DAEs in which the algebraic constraint is treated by the Rosenbrock partition.  
\ifreport
To better illustrate the structure of the method, we provide the stage computations below:
\begin{align*}
    \mathbf{k}^{\{1+2\}}_1 &= h \, \fun^{\{1\}}(\y_n) + h \, \fun^{\{2\}}(\y_n), \\
    \mathbf{k}^{\{3\}}_1 &= h \, \fun^{\{3\}}(\y_n) + h \, \gamma \, \Jac_n^{\{3\}} \left( \mathbf{k}^{\{1+2\}}_1 + \mathbf{k}^{\{3\}}_1 \right), \\
    \mathbf{k}^{\{1+2\}}_2 &= h \, \fun^{\{1\}} \mleft( \y_{n} + \mathbf{k}^{\{1+2\}}_1 + \mathbf{k}^{\{3\}}_1 \mright) + h \, \fun^{\{2\}} \mleft( \y_{n} + \frac{1}{2} \left( \mathbf{k}^{\{1+2\}}_1  + \mathbf{k}^{\{1+2\}}_2 \right) + \mathbf{k}^{\{3\}}_1 \mright), \\
    \mathbf{k}^{\{3\}}_2 &= h \, \fun^{\{3\}} \mleft( \y_{n} + \mathbf{k}^{\{1+2\}}_1 + \mathbf{k}^{\{3\}}_1 \mright) + h \, \gamma \, \Jac_n^{\{3\}} \left( \mathbf{k}^{\{1+2\}}_2 - \mathbf{k}^{\{1+2\}}_1 + \mathbf{k}^{\{3\}}_2 - \mathbf{k}^{\{3\}}_1 \right), \\
    \y_{n+1} &= \y_n + \frac{1}{2} \left( \mathbf{k}^{\{1+2\}}_1  + \mathbf{k}^{\{1+2\}}_2 \right) + (1 - \gamma) \, \mathbf{k}^{\{3\}}_1 + \gamma \, \mathbf{k}^{\{3\}}_2.
\end{align*}
\fi
The implicit and explicit trapezoidal rules share the same $\b$, which allows us to use the combined stage $\mathbf{k}^{\{1+2\}}_i = \mathbf{k}^{\{1\}}_i + \mathbf{k}^{\{2\}}_i$ as discussed in \cref{rem:imex-same-b}.  Note that when $\fun^{\{2\}}(\y) = 0$, the method degenerates into a two-way partitioned IMEX GARK-ROS scheme which we refer to as IMEX-ROS22.

\subsection{Third order IMEX GARK-ROW schemes}

We explore IMEX GARK-Rosenbrock-W methods that are suitable for index-1 DAEs and are equipped with an embedded method for error estimation and control.  The special cases described in \cref{rem:imex-special-case,rem:imex-same-b} are used to reduce the number of coefficients and order conditions.

We first consider the case when $s\E = s\I = 4$.  For the base Rosenbrock method, we enforce traditional ROW and DAE order conditions up to order three.  Similarly, the explicit base method must satisfy Runge--Kutta order conditions up to order three.  These base methods share the embedded coefficients $\widehat{\b}$, which must give a solution of order two.  To form an IMEX pair, the coupling condition \cref{eqn:GARK-ROW-IMEX-case1-order3} and DAE coupling condition \cref{eqn:GARK-ROS-IMEX-DAE-case1-order3} are imposed.  There are still several free parameters left after solving these order conditions, and in our method derivation procedure, they are used to optimize the stability and principal error.  Our method, IMEX-ROW3(2)4, pairs the explicit Runge--Kutta scheme
\begin{subequations}
\begin{equation} \label{eqn:IMEX-ROW3(2)4:RK}
\scalebox{0.85}{$
        \begin{butchertableau}{c|cccc}
             {\scriptstyle 0} & {\scriptstyle 0} & {\scriptstyle 0} & {\scriptstyle 0} & {\scriptstyle 0} \\
             {\scriptstyle 2 \gamma}  & {\scriptstyle 2 \gamma}  & {\scriptstyle 0} & {\scriptstyle 0} & {\scriptstyle 0} \\
             \frac{\gamma + 1}{2} & -\frac{15 \gamma ^2}{16}+\frac{103 \gamma }{32}-\frac{5}{8} & \frac{15 \gamma ^2}{16}-\frac{87 \gamma }{32}+\frac{9}{8} & {\scriptstyle 0} & {\scriptstyle 0} \\
             {\scriptstyle 1} & -\frac{81 \gamma ^2}{272}+\frac{111 \gamma }{136}+\frac{265}{544} & \frac{\gamma ^2}{16}+\frac{\gamma }{8}-\frac{25}{32} & \frac{4 \gamma ^2}{17}-\frac{16 \gamma }{17}+\frac{22}{17} & {\scriptstyle 0} \\ \hline
             & -\frac{9 \gamma ^2}{34}+\frac{19 \gamma }{34}+\frac{3}{68} & \frac{5 \gamma ^2}{2}-\frac{13 \gamma }{2}+\frac{5}{4} & -\frac{38 \gamma ^2}{17}+\frac{84 \gamma }{17}-\frac{5}{17} & {\scriptstyle \gamma} \\ \hline
             & -\frac{57 \gamma ^2}{272}+\frac{109 \gamma }{272}+\frac{9}{136} & \frac{47 \gamma ^2}{16}-\frac{31 \gamma }{4}+\frac{23}{16} & -\frac{40 \gamma ^2}{17}+\frac{201 \gamma }{34}-\frac{15}{34} & -\frac{3 \gamma ^2}{8}+\frac{23 \gamma }{16}-\frac{1}{16} \\
        \end{butchertableau}
$}
\end{equation}
with an L-stable Rosenbrock-W method with coefficients
\begin{equation} \label{eqn:IMEX-ROW3(2)4:ROW}
    \begin{split}
&\scalebox{0.85}{$
    \boldsymbol{\alpha} = \begin{butchermatrix}
         0 & {\scriptstyle 0} & {\scriptstyle 0} & {\scriptstyle 0} \\
         {\scriptstyle 2 \gamma}  & {\scriptstyle 0} & {\scriptstyle 0} & {\scriptstyle 0} \\
         -\frac{9 \gamma ^2}{8}+\frac{115 \gamma }{32}-\frac{19}{32} & \frac{9 \gamma ^2}{8}-\frac{99 \gamma }{32}+\frac{35}{32} & {\scriptstyle 0} & {\scriptstyle 0} \\
         \frac{9 \gamma ^2}{34}-\frac{19 \gamma }{34}+\frac{31}{68} & -\frac{\gamma ^2}{2}+\frac{3 \gamma }{2}-\frac{3}{4} & \frac{4 \gamma ^2}{17}-\frac{16 \gamma }{17}+\frac{22}{17} & {\scriptstyle 0} \\
    \end{butchermatrix}, 
    $}
    \\
&\scalebox{0.85}{$
    \boldsymbol{\gamma} = \begin{butchermatrix}
         \gamma  & {\scriptstyle 0} & {\scriptstyle 0} & {\scriptstyle 0} \\
         -2 \gamma  & \gamma  & {\scriptstyle 0} & {\scriptstyle 0} \\
         \frac{3 \gamma ^2}{2}-\frac{157 \gamma }{32}+\frac{33}{32} & -\frac{3 \gamma ^2}{4}+\frac{57 \gamma }{32}-\frac{21}{32} & \gamma  & {\scriptstyle 0} \\
         -\frac{9 \gamma ^2}{17}+\frac{19 \gamma }{17}-\frac{7}{17} & 3 \gamma ^2-8 \gamma +2 & -\frac{42 \gamma ^2}{17}+\frac{100 \gamma }{17}-\frac{27}{17} & \gamma  \\
    \end{butchermatrix},
    $}
    \end{split}
\end{equation}
\end{subequations}
where $\gamma \approx 0.44$ is the middle root of $6 \gamma ^3-18 \gamma ^2+9 \gamma -1 = 0$.  The $\b$ and $\widehat{\b}$ coefficients in \cref{eqn:IMEX-ROW3(2)4:ROW} are the same as in \cref{eqn:IMEX-ROW3(2)4:RK}.  Thanks to the stiff accuracy of the Rosenbrock method, \cref{eqn:inconsistent-ic-z-1} is satisfied as well; however, we were unable to cancel higher order error terms for inconsistent initial conditions.

We also derive a third order scheme with $s\E = s\I = 5$ as it affords a smaller $\gamma_{i,i}$ and sufficient degrees of freedom to satisfy \cref{eqn:inconsistent-ic-z-2,eqn:inconsistent-ic-x-1}, thus eliminating errors associated with inconsistent initial values up to $\mathcal{O}(h\delta)$.  On top of the simplifying assumptions and order conditions used with four stages, we take $\boldsymbol{\alpha}\E = \boldsymbol{\alpha}\I$, such that the method looks like an unpartitioned Rosenbrock-W method with $\Lb = \fun\I_\y$.   For DAEs however, one cannot expect a general Rosenbrock-W method to attain full order when the Jacobian of $\fun\Alg$ is used; the order condition \cref{eqn:GARK-ROS-IMEX-DAE-case1-order3} is required for this.

Based on the aforementioned constraints, our five-stage method, named IMEX-ROW3(2)5, has the coefficients
\begin{equation} \label{eqn:IMEX-ROW3(2)5:ROW}
    \begin{split}
\scalebox{0.85}{$
    \boldsymbol{\alpha}\E = \boldsymbol{\alpha}\I = \begin{butchermatrix}
         {\scriptstyle 0} & {\scriptstyle 0} & {\scriptstyle 0} & {\scriptstyle 0} & {\scriptstyle 0} \\
         \frac{1}{2} & {\scriptstyle 0} & {\scriptstyle 0} & {\scriptstyle 0} & {\scriptstyle 0} \\
         \frac{5062}{13725} & \frac{4088}{13725} & {\scriptstyle 0} & {\scriptstyle 0} & {\scriptstyle 0} \\
         \frac{173067}{636265} & \frac{495828}{636265} & -\frac{24705}{127253} & {\scriptstyle 0} & {\scriptstyle 0} \\
         \frac{30859}{262800} & -\frac{547}{21900} & \frac{183}{146} & -\frac{18179}{52560} & {\scriptstyle 0}
    \end{butchermatrix}, 
    \qquad
    \b = \begin{butchermatrix}
        \frac{5225}{21024} \\ -\frac{407}{2190} \\ \frac{6039}{4672} \\ -\frac{127253}{210240} \\ \frac{1}{4}
    \end{butchermatrix},
    $} \\
\scalebox{0.85}{$
    \boldsymbol{\gamma} = \begin{butchermatrix}
         \frac{1}{4} & {\scriptstyle 0} & {\scriptstyle 0} & {\scriptstyle 0} & {\scriptstyle 0} \\
         -\frac{1}{2} & \frac{1}{4} & {\scriptstyle 0} & {\scriptstyle 0} & {\scriptstyle 0} \\
         -\frac{4762}{13725} & -\frac{2563}{13725} & \frac{1}{4} & {\scriptstyle 0} & {\scriptstyle 0} \\
         -\frac{156792}{636265} & -\frac{685353}{636265} & \frac{82350}{127253} & \frac{1}{4} & {\scriptstyle 0} \\
         \frac{22969}{175200} & -\frac{3523}{21900} & \frac{183}{4672} & -\frac{18179}{70080} & \frac{1}{4}
    \end{butchermatrix},  \qquad
    \widehat{\b} = \begin{butchermatrix}
        \frac{9095}{539616} \\  \frac{27387}{56210}  \\  \frac{421083}{359744}  \\ -\frac{812861}{770880}  \\  \frac{117}{308}
    \end{butchermatrix}.
    $}
    %
    \end{split}
\end{equation}
When viewed as an unpartitioned Rosen\-brock-W method, IMEX-ROW3(2)5 is stiffly accurate and L-stable.

\subsection{Fourth order IMEX GARK-ROS scheme}

Order four introduces significantly more order conditions, and it appears six stages is the minimum required for an IMEX GARK-ROS scheme that is suitable for index-1 DAEs and includes an embedded method.  For the base ROS method, classical and DAE order conditions up to order four are necessary, but we include ROW order conditions up to order three as well.  The base Runge--Kutta method uses Butcher's first column simplifying assumption $D(1)$ \citep{Butcher_1964_IRK}, which leaves fives order conditions to achieve order four.  With \cref{rem:imex-special-case,rem:imex-same-b}, the IMEX coupling conditions are \cref{eqn:GARK-ROW-IMEX-case1-order3,eqn:GARK-ROW-IMEX-case1-order4}, and the DAE coupling conditions are \cref{eqn:GARK-ROS-IMEX-DAE-case1-order3,eqn:GARK-ROS-IMEX-DAE-case1-order4}.  The embedded method, with coefficients $\widehat{\b}$, must satisfy all these order conditions to one order lower.  We solve the order conditions and use remaining free coefficients for tuning stability and principal error.  
The final method, IMEX-ROS4(3)6, pairs the explicit Runge--Kutta scheme
\begin{equation*}
\scalebox{0.85}{$
    \begin{butchertableau}{c|cccccc}
        {\scriptstyle 0} & {\scriptstyle 0} & {\scriptstyle 0} & {\scriptstyle 0} & {\scriptstyle 0} & {\scriptstyle 0} & {\scriptstyle 0} \\
        \frac{1}{2} & \frac{1}{2} & {\scriptstyle 0} & {\scriptstyle 0} & {\scriptstyle 0} & {\scriptstyle 0} & {\scriptstyle 0} \\
        \frac{9}{10} & \frac{4761}{11050} & \frac{2592}{5525} & {\scriptstyle 0} & {\scriptstyle 0} & {\scriptstyle 0} & {\scriptstyle 0} \\
        \frac{2}{5} & \frac{3779}{99450} & \frac{12931}{44200} & \frac{5}{72} & {\scriptstyle 0} & {\scriptstyle 0} & {\scriptstyle 0} \\
        \frac{5}{6} & -\frac{9468553}{45647550} & \frac{18193697}{30431700} & -\frac{92843}{413100} & \frac{1352}{2025} & {\scriptstyle 0} & {\scriptstyle 0} \\
        {\scriptstyle 1} & \frac{5613193}{5967000} & \frac{261179}{884000} & \frac{18091}{108000} & -\frac{13609}{19500} & \frac{153}{520} & {\scriptstyle 0} \\ \hline
        & \frac{113}{720} & \frac{37}{96} & -\frac{125}{288} & \frac{125}{624} & \frac{459}{1040} & \frac{1}{4} \\ \hline
        & \frac{433321}{3204900} & \frac{121913}{569760} & -\frac{25667}{1025568} & \frac{6024}{15431} & \frac{965889}{6172400} & \frac{1531}{11870}
    \end{butchertableau}
    $}
\end{equation*}
with the stiffly accurate, L-stable Rosenbrock scheme with coefficients
\begin{align*}
&\scalebox{0.85}{$
    \boldsymbol{\alpha} = \begin{butchermatrix}
        {\scriptstyle 0} & {\scriptstyle 0} & {\scriptstyle 0} & {\scriptstyle 0} & {\scriptstyle 0} & {\scriptstyle 0} \\
        \frac{1}{2} & {\scriptstyle 0} & {\scriptstyle 0} & {\scriptstyle 0} & {\scriptstyle 0} & {\scriptstyle 0} \\
        \frac{87}{140} & \frac{39}{140} & {\scriptstyle 0} & {\scriptstyle 0} & {\scriptstyle 0} & {\scriptstyle 0} \\
        -\frac{331}{1260} & \frac{17}{28} & \frac{1}{18} & {\scriptstyle 0} & {\scriptstyle 0} & {\scriptstyle 0} \\
        \frac{84025}{231336} & -\frac{755}{9639} & -\frac{425}{1944} & \frac{4225}{5508} & {\scriptstyle 0} & {\scriptstyle 0} \\
        \frac{1091}{2160} & \frac{29}{32} & \frac{145}{864} & -\frac{545}{624} & \frac{153}{520} & {\scriptstyle 0}
    \end{butchermatrix}, 
    $}
    \quad
\scalebox{0.85}{$
    \boldsymbol{\gamma} = \begin{butchermatrix}
        \frac{1}{4} & {\scriptstyle 0} & {\scriptstyle 0} & {\scriptstyle 0} & {\scriptstyle 0} & {\scriptstyle 0} \\
        -\frac{1}{2} & \frac{1}{4} & {\scriptstyle 0} & {\scriptstyle 0} & {\scriptstyle 0} & {\scriptstyle 0} \\
        -\frac{183}{700} & \frac{57}{700} & \frac{1}{4} & {\scriptstyle 0} & {\scriptstyle 0} & {\scriptstyle 0} \\
        \frac{257}{700} & -\frac{731}{1400} & -\frac{1}{8} & \frac{1}{4} & {\scriptstyle 0} & {\scriptstyle 0} \\
        \frac{33925}{231336} & \frac{45835}{77112} & \frac{2725}{16524} & -\frac{1300}{1377} & \frac{1}{4} & {\scriptstyle 0} \\
        -\frac{47}{135} & -\frac{25}{48} & -\frac{65}{108} & \frac{335}{312} & \frac{153}{1040} & \frac{1}{4}
    \end{butchermatrix}. $}
\end{align*}



\section{Numerical Experiments}
\label{sec:experiments}

In this section, we present the results from two numerical experiments that verify the linearly-implicit GARK order condition theory and the convergence properties of the methods derived in \cref{sec:methods}.  For comparison, we also test Ros-ERK3,2 from \citep{Hai_2014_ERK-ROS}, ROS3P from \citep{Lang2001}, and RODAS from \citep[Section VI.4]{Hairer_book_II}.

\subsection{Brusselator reaction-diffusion PDE}

The problem BRUSS from \citep[pg 148]{Hairer_book_II}, is a one-dimensional reaction-diffusion problem governed by the equations
\begin{equation} \label{eqn:bruss}
    \begin{split}
        \frac{\partial u}{\partial t} &= A + u^2 \, v - (B + 1) \, u + \alpha \, \frac{\partial^2 u}{\partial x^2}, \qquad
        \frac{\partial v}{\partial t} = B \, u - u^2 \, v + \alpha \, \frac{\partial^2 v}{\partial x^2},
    \end{split}
\end{equation}
with $A = 1$, $B = 3$, and $\alpha = 1 / 50$.  The spatial domain is $x \in [0, 1]$, and the time domain is $t \in [0, 10]$. The boundary and initial conditions are
\begin{align*}
    &u(x=0, t) = u(x=1, t) = 1,
    &&v(x=0, t) = v(x=1, t) = 3; \\
    &u(x, t= 0) = 1 + \sin(2 \, \pi \, x),
    &&v(x, t=0) = 3.
\end{align*}
Second order central finite differences are applied to discretize the spatial dimension on a uniform grid with $N = 500$ interior points.

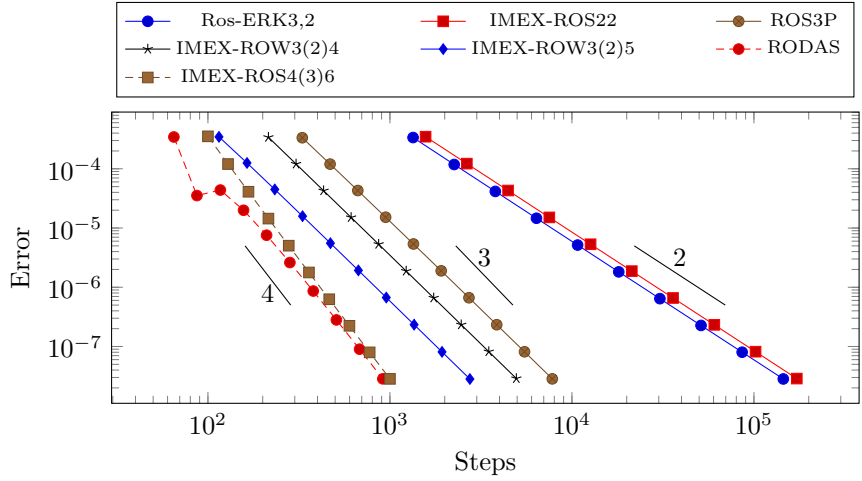
\begin{figure}[ht!]
    \centering
    \begin{tikzpicture}
        \begin{loglogaxis}[
    	        legend entries={{Ros-ERK3,2},IMEX-ROS22,ROS3P,IMEX-ROW3(2)4,IMEX-ROW3(2)5,RODAS,IMEX-ROS4(3)6},
        	    xlabel=Steps,
        	    ylabel=Error,
        	    width=0.95\linewidth,
        	    height=0.45\linewidth
            ]
            \addplot table[x index=0,y index=7] {bruss.dat};
            \addplot table[x index=1,y index=8] {bruss.dat};
            \draw (axis cs:2.2e4,5e-6) -- node[above]{$2$} (axis cs:6.96e4,5e-7);
            \addplot table[x index=2,y index=9] {bruss.dat};
            \addplot table[x index=3,y index=10] {bruss.dat};
            \addplot table[x index=4,y index=11] {bruss.dat};
            \draw (axis cs:2.3e3,5e-6) -- node[above]{$3$} (axis cs:4.74e3,5e-7);
            \addplot table[x index=5,y index=12] {bruss.dat};
            \addplot table[x index=6,y index=13] {bruss.dat};
            \draw (axis cs:1.6e2,5e-6) -- node[below]{$4$} (axis cs:2.85e2,5e-7);
        \end{loglogaxis}
    \end{tikzpicture}
    \caption{IMEX convergence results on the Brusselator problem \cref{eqn:bruss}.}
    \label{fig:bruss-convergence}
\end{figure}

The stiffness in \cref{eqn:bruss} primarily comes from the diffusion terms.  Therefore, we treat them linearly implicitly and the remaining reaction terms explicitly.  For each of the itegrators, we compute the numerical error for a range of ten step sizes.  Error is measured as the two-norm of the difference of the numerical solution and a highly accurate reference solution at $t = 10$.  The converge plots are shown in \cref{fig:bruss-convergence}.  In all cases, the numerical orders of convergence match the theoretical ones.

\subsection{ZLA-kinetics problem}

The ZLA-kinetics problem is a nonlinear index-1 DAE modelling the reaction of two chemicals as carbon dioxide is added to the system.  A detailed description of this problem and its origin is provided in \citep{stortelder1998parameter}.  It is governed by the following five differential equations and one algebraic constraint:
\begin{equation} 
\label{eqn:zla}
    \begin{alignedat}{2}
        y_1' &= -2 \, r_1 + r_2 - r_3 - r_4 , & \qquad
        y_2' &= -\frac{1}{2} \, r_1 - r_4 - \frac{1}{2} \, r_5 + F_{\rm in}, \\
        y_3' &= r_1 - r_2 + r_3, & \qquad
        y_4' &= -r_2 + r_3 -2 \, r_4, \\
        y_5' &= r_2 - r_3 + r_5, & \qquad
        0 &= K_s \, y_1 \, y_4 - y_6.
    \end{alignedat}
\end{equation}
The auxiliary variables and parameters are defined as:
\begin{alignat*}{5}
    r_1 &= k_1 \, y_1^4 \, y_2^{1/2},
    \quad & r_2 &= k_2 \, y_3 \, y_4,
    \quad & r_3 &= (k_2/K) \, y_1 \, y_5, \\
    r_4 &= k_3 \, y_1 \, y_4^2,
    \quad & r_5 &= k_4 \, y_6^2 \, y_2^{1/2},
    \quad & F_{\rm in} &= klA \left( p(\text{CO}_2)/H - y_2 \right), \\
%
    k_1 &= 18.7,
    & \quad k_2 &= 0.58,
    & \quad k_3 &= 0.09, \\
    k_4 &= 0.42,
    & \quad K &= 34.4,
    & \quad klA &= 3.3, \\
    K_s &= 115.83,
    & \quad p(\text{CO}_2) &= 0.9,
    & \quad H &= 737.
\end{alignat*}
The system is integrated from $t = 0$ to $t = 180$ starting from the initial value
\begin{equation*}
    \y(t=0) = \begin{bmatrix}
        0.444 & 0.00123 & 0 & 0.007 & 0 & K_s \, y_{0,1} \, y_{0,4}
    \end{bmatrix}^T,
\end{equation*}
which is consistent with the algebraic constraint.

\begin{figure}[ht!]
    \centering
    \begin{tikzpicture}
        \begin{loglogaxis}[
    	        legend entries={{Ros-ERK3,2},IMEX-ROS22,ROS3P,IMEX-ROW3(2)4,IMEX-ROW3(2)5,RODAS,IMEX-ROS4(3)6},
        	    xlabel=Steps,
        	    ylabel=Error,
        	    width=0.95\linewidth,
        	    height=0.45\linewidth
            ]
            \addplot table[x index=0,y index=7] {zla.dat};
            \addplot table[x index=1,y index=8] {zla.dat};
            \draw (axis cs:9e5,1e-11) -- node[above]{$2$} (axis cs:28.46e5,1e-12);
            \addplot table[x index=2,y index=9] {zla.dat};
            \addplot table[x index=3,y index=10] {zla.dat};
            \addplot table[x index=4,y index=11] {zla.dat};
            \draw (axis cs:4e4,1e-11) -- node[above]{$3$} (axis cs:8.62e4,1e-12);
            \addplot table[x index=5,y index=12] {zla.dat};
            \addplot table[x index=6,y index=13] {zla.dat};
            \draw (axis cs:7e3,1e-11) -- node[below]{$4$} (axis cs:1.24e4,1e-12);
        \end{loglogaxis}
    \end{tikzpicture}
    \caption{IMEX convergence results on the ZLA-kinetics problem \cref{eqn:zla}.}
    \label{fig:ZLA-convergence}
\end{figure}
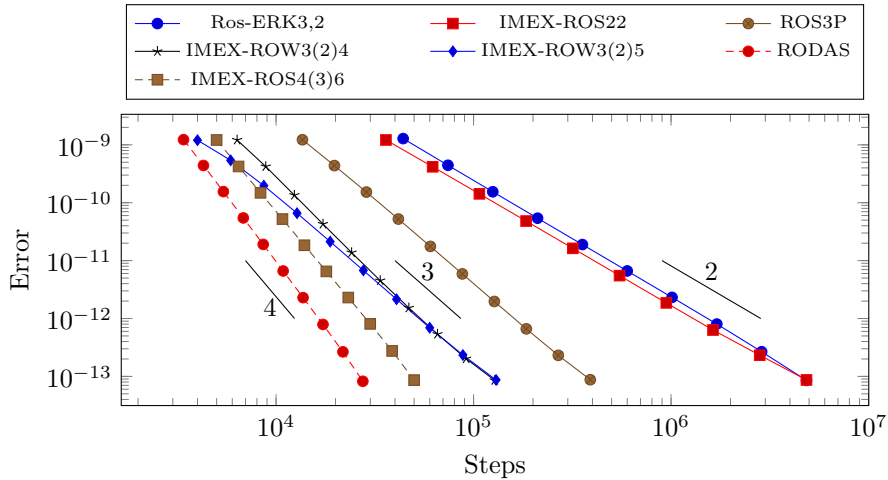

We use the ZLA-kinetics problem to verify DAE convergence properties of the IMEX methods proposed in \cref{sec:methods}.  In the numerical experiment, the differential variables are treated explicitly, while the algebraic variable is treated linearly implicitly.  \Cref{fig:ZLA-convergence} plots the error versus the number of steps taken to solve the DAE.  Like the Brusselator experiment, error is measured in the two-norm with respect to a reference solution.  All methods achieve their theoretical orders of convergence.

\section{Discussion}
\label{sec:discussion}

This paper constructs new families of linearly implicit multimethods. The authors' GARK framework extends traditional Runge--Kutta schemes to multimethods suitable for the discretization of multiphysics systems. In a similar vein, the GARK-ROS/GARK-ROW framework extends traditional Ro\-sen\-brock/Ro\-sen\-brock-W schemes to multimethods. 

A general order conditions theory for linearly implicit methods with any number of partitions, using exact or approximate Jacobians, is developed using B-series over the sets of $\NT$ trees (for exact Jacobian) and $\NTW$ trees  (for inexact Jacobians). Order conditions for the solution of two-way partitioned index-1 differential-algebraic equations are developed using B-series over the set of $\DAT$ trees. We use the framework to develop decoupled linearly implicit schemes, which treat implicitly one process at a time; linearly implicit/explicit methods, which treat one process explicitly and one implicitly; and linearly implicit/explicit/implicit methods that discretize some processes with Rosenbrock schemes, other with diagonally implicit Runge--Kutta schemes, and other with explicit Runge--Kutta schemes. Practical GARK-ROS and GARK-ROW schemes of orders two, three, and four are constructed.



\bibliographystyle{plain}


 
\end{document}